\documentclass[a4paper,12pt]{article}
\author{G\'erard Ben Arous\thanks{Courant Institute of Mathematical Sciences, New York University. Research supported in part by the National Science Foundation under grants DMS-0806180 and OISE-0730136.}, Alan Hammond \thanks{Department of Statistics, Oxford University. Research supported in part by NSF grants DMS-0806180 and OISE-0730136 and by EPSRC grant EP/I004378/1. This research was initiated while the second author was at N.Y.U.}}
\title{Randomly biased walks on subcritical trees}

\usepackage{amsfonts,psfrag,epsfig}
\usepackage{appendix}
\usepackage{url}

\def\P{\mathbb{P}}
\def\Z{\mathbb{Z}}
\def\N{\mathbb{N}}
\def\E{\mathbb{E}}

\def\R{\mathbb{R}}
\def\binf{q}
\def\bsup{Q}

\def\coup{\Theta}

\def\cpo{d_1}

\def\cpob{d_2}

\def\barecon{B}
\def\conec{D_1}

\def\trwtbd{v_0}
\def\wgt{weighted }

\def\coup{\mathbb{Q}}
\def\basedef{{\rm base}}
\def\vmax{v_{\rm max}}
\def\vlast{{v_{\rm base}}}
\def\ancestor{{v_{\rm anc}}}
\def\fbp{{v_{\rm fbp}}}
\def\vnec{v_{\rm ess}}
\def\vchild{{v_{\rm child}}}

\def\ef{EF}
\def\opf{OF}
\def\phni{\P_{h,\nu}^{(i)}}
\def\phniplus{\P_{h,\nu}^{(i),+}}
\def\fso{{\rm FSO}}
\def\ess{{\rm ESS}}
\def\vbase{{v_{\rm base}}}

\def\conlem{c}

\def\conbare{B}
\def\constd{d}
\def\falldeep{\mathcal{DE}}
\def\pfd{p_{\rm de}}
\def\omegachild{\omega_*}
\def\omegachildbig{\omega_*}
\def\czerosm{c_1}
\def\conesm{c_2}
\def\ctwosm{c_3}
\def\cone{c_4}
\def\conom{C_1}
\def\ctwo{C_2}
\def\conthrbig{C_3}
\def\reo{E_0}
\def\reotilde{E}
\def\ezero{\mathcal{E}_1}
\def\eone{\mathcal{E}_2}
\def\errzero{\mathcal{E}_0}
\newtheorem{definition}{Definition}
\newtheorem{lemma}{Lemma}
\newtheorem{corollary}{Corollary}
\newtheorem{theorem}{Theorem}
\newtheorem{prop}{Proposition}

\newtheorem{hyp}{Hypothesis}

\def\build#1_#2^#3{\mathrel{ \mathop{\kern 0pt#1}\limits_{#2}^{#3}}}

\begin{document}
%\nocite{*}
\maketitle

\begin{abstract}
As a model of trapping by biased motion in random structure, 
we study the time taken for a biased random walk to return to the root of a subcritical Galton-Watson tree. 
We do so for trees in which these biases are randomly chosen, independently for distinct edges, according to a law 
that satisfies a logarithmic non-lattice condition. The mean return time of the walk is in essence given by  
the total conductance of the tree. We determine the asymptotic decay of this total conductance, finding it to have 
a pure power-law decay. In the case of the conductance associated to a single vertex at maximal depth in the tree, 
 this asymptotic decay may be analysed by the classical defective renewal theorem, due to the non-lattice edge-bias assumption.   
  However, the derivation of the decay for total conductance requires computing an additional constant multiple outside the power-law that allows for the contribution of all vertices close to the base of the tree. This computation entails a detailed study of a convenient decomposition of the tree, under conditioning on the tree having high total conductance.
As such,
our principal conclusion may be viewed as a development of renewal theory in the context of random environments. 

For randomly biased random walk on a supercritical Galton-Watson tree with positive extinction probability, 
our main results may be regarded as a description of the slowdown mechanism caused by the presence of subcritical trees adjacent to the backbone that may act as traps that detain the walker. Indeed, this conclusion is exploited in \cite{GerardAlan} to obtain a stable limiting law for walker displacement in such a tree. 
%In these terms, it is the role of the non-lattice condition on edge-biases to eliminate a log-periodic effect in the total conductance of traps that would otherwise obtain. 
\end{abstract}

%\small{
%\tableofcontents
%}
%\normalsize

\begin{section}{Introduction}

The inquiry into drift and trapping for  biased random walks on random structures
has been pursued both by physicists and mathematicians. 
 It was noted long ago in the physics literature \cite{BB} that, in disordered media, due to trapping in dead-end branches, the mean velocity would not be a monotone function of the bias. It was then pointed out that, in fact,  strong biases could even produce sub-ballistic regimes, i.e., with zero velocity; (see \cite{Dhar}, \cite{DS}, or \cite{Havlin} for a physics survey of the general phenomena of anomalous diffusion). This hypothesis was confirmed for biased random walks on random trees by Lyons, Pemantle and Peres \cite{LPP96}, and on supercritical percolation clusters, independently by Berger, Gantert and Peres in dimension two \cite{BGP}, and by Sznitman in all dimensions $d \geq 2$ \cite{Sznitman1}. For $d \geq 2$, the existence of a critical bias separating the sub-ballistic and ballistic regimes will be demonstrated in \cite{AlexAlan}. 
These studies naturally raise the question of understanding more deeply the means by which the walk is slowed down  by trapping structures in this sub-ballistic regime. A. Sznitman mentions in \cite{Sznitman2} that this mechanism  seems similar to that responsible for aging in Bouchaud's trap model (\cite{BAC},\cite{Bouchaud}) whose main features are that the mean trapping times have a power-law distribution and that the tail of the distribution of the trapping times, conditionally on their mean, is exponential.

For random walks on supercritical Galton-Watson trees whose bias away from the root is constant, 
such a study of the trapping phenomenon was undertaken  in \cite{BAFGH}. In this context, the dead-ends responsible for the slowing down of the walk are simply subcritical Galton-Watson trees that hang off the backbone of the supercritical tree. Thus, the trapping times are related to return times to the root for biased random walks on subcritical trees. The paper \cite{BAFGH} analyses a log-periodicity phenomenon which is very reminiscent of the classical lattice effect for random walks in random environments \cite{KKS}, \cite{ESZ1}, \cite{ESZ2}, \cite{ZeitouniStFlour}.  This lattice effect, which has been discussed in the physics literature \cite{StaSto}, causes holding times in traps to cluster around powers of the bias parameter. That is to say, 
although it may seem reasonable that such a walk has a scaling limit similar to the one obtained for Bouchaud trap models that is discussed in \cite{BAC} and \cite{BACerny}, the lattice effect is a discrete inhomogeniety in walker displacement that is persistent on all time scales and that prevents the existence of such a scaling limit.

In this paper, we study the question of how long it takes a biased random walk in a subcritical Galton-Watson tree to return to the root of the tree. In the context of biased walks on supercritical trees, our inquiry amounts to an investigation of the nature of the delaying mechanism of the subcritical trees that act as potential traps for the walk. We undertake the inquiry in the case that the biases on edges that determine the distribution of the walk are random, rather than identically equal, imposing a non-lattice condition on this randomness that serves to eliminate the log-periodic effect from the return-time distribution. 

This paper is the first of two. Indeed, its sequel \cite{GerardAlan} exploits
the understanding of the trapping mechanism that we obtain in the present article to prove a stable scaling limit for the randomly biased random walk on supercritical Galton-Watson trees with leaves; in this way, the two papers 
form a counterpart to \cite{BAFGH}, in that they show how edge-bias randomization dissipates the persistent discrete inhomogeneity that obtains in the case of constant bias.

Randomly biased random walks on infinite trees have been studied by \cite{LyonsPemantle}, 
who provide a criterion for recurrence, or transience, valid for a very broad class of trees. 
 In the case of supercritical Galton-Watson trees, randomly biased walk is studied in the recurrent regime by \cite{FHS}, who investigate the high-$n$ asymptotic of the maximal displacement of the walk during $[0,n]$. In \cite{aidekon}, the same model is analysed in the case of zero extinction probability and in the transient regime, the author presenting criteria for zero speed and for positive speed. The tree being leafless in \cite{aidekon}, this work is not concerned with trapping, but rather is a contribution to understanding how increasing bias causes more rapid progress for the walk on the backbone of a supercritical Galton-Watson tree. 
Regarding this last question, the monotonicity of speed as a function of bias is not known even when the bias is a constant: see Question 2.1 of \cite{LPP97}. %\cite{LPP95}

\bigskip

We now define the general class of biased random walks on finite rooted trees in which we are interested.
Let $T$ be a finite rooted tree, with vertex set $V(T)$, edge set $E(T)$, and root $\phi$. Such a tree will be called weighted if, to each unoriented edge $e \in E(T)$  is assigned a number $\beta_e \in [\binf,\bsup]$, that we will call the bias of the edge $e$. Here, $\bsup \geq \binf > 1$ are fixed constants. In this paper, in referring to a tree, we will generally mean a weighted tree. We  define  the law $\P_{T,\beta}$ of the $\beta$-biased random walk $\big\{ X_i: i \in \N \big\}$ on the set of vertices of the tree $T$ as a Markov chain on $V(T)$ with the following transition rules.

\begin{definition}\label{defwalk}
If a vertex $v \in V(T)$, $v
\not= \phi$, has offspring $v_1,\ldots,v_k$, then, for each $n \in \N$,  
\begin{eqnarray}
\P_{T,\beta} \Big( X_{n + 1} = \overleftarrow{v} \Big\vert X_n = v \Big) & = &
\frac{1}{1 + \sum_{i=1}^k \beta_{v,v_i}},
        \nonumber \\
\P_{T,\beta} \Big( X_{n + 1} = v_j \Big\vert X_n = v \Big)    
& = & \frac{\beta_{v,v_j}}{1 + \sum_{i=1}^k \beta_{v,v_i}},
 \qquad \textrm{for $1 \leq j \leq k$}, \nonumber
\end{eqnarray}
where $\overleftarrow{v}$ denotes the parent of $v$, i.e., the neighbor of $v$ which is the closest to the root. The jump-law from the root is given by
$$
\P_{T,\beta} \Big( X_{n + 1} = \phi_j \Big\vert X_n = \phi \Big)    
 =  \frac{\beta_{\phi,\phi_j}}{\sum_{i=1}^k \beta_{\phi,\phi_i}},
 \qquad \textrm{for $1 \leq j \leq k(\phi)$}.
$$  
For $v \in V(T)$, we write $\P_{T,\beta}^v$ for the law of $\P_{T,\beta}$ given that
$X(0) = v$. 
We write $\E_{T,\beta}^v$ for the expectation value under $\P_{T,\beta}^v$. 
\end{definition}

It is easy to see these biased random walks are reversible and to compute their invariant measures. 
Indeed, they fall in the general class of ``random walks on weighted graphs''(see \cite{Kumagai}, and \cite{Barlow} for  recent expositions).
\begin{definition}\label{defomega}
For every vertex $v \in V(T)$, define $P_{\phi,v}$ to be the unique simple path from the root to $v$, and $\omega(v)$ to be the product of the biases $\beta_e$ of the edges $e$
along the path $P_{\phi,v}$.  If   $e=(v,w)$ is an edge of $T$, with $v=\overleftarrow{w}$, we define the weight (or conductance) of $e$ by $\mu(e)= \mu(v,w)= \omega(w)$ . Finally, we define the unnormalized measure $\mu$ on $V(T)$ by setting $\mu(v)$ equal to the sum of the conductances of the edges adjacent to the vertex $v \in V(T)$.
 \end{definition}

It is clear that the measure $\mu$ is reversible and thus invariant for the Markov chain $X_n$.
Moreover, by the general random walks on weighted graphs, the transition rules between the adjacent vertices $v$ and $w$ can simply be rewritten in terms of the  conductance of the edge between $v$ and $w$ as well as the measure $\mu(v)$:

\begin{equation}
\P_{T,\beta} \Big( X_{n + 1} = w \Big\vert X_n = v \Big)  = 
\frac{\mu(v,w)}{ \mu(v)}.
\end{equation}

We will be interested in the first return time to the root of the walk that starts at the root,  where the return time is naturally defined as $\inf \big\{ n >1 , X_n=\phi \big\}$ under the law $\P_{T,\beta}^\phi$.

As we have mentioned, the random trees that we consider are subcritical Galton-Watson trees. 
\begin{definition}\label{deftreelaw}
Let $h = \big\{ h_i: i \in \N \big\}$, $\sum_{i=0}^\infty h_i = 1$, denote an offspring distribution. 
We write $\P_h$ for the Galton-Watson tree with offspring distribution $h$, that is, for the law on rooted trees in which each vertex has an independent and $h$-distributed number of offspring.
\end{definition}

We will make the following assumption on the offspring distribution.

\begin{hyp}\label{hyph}
Let $\big\{ h_i:  i \in \N \big\}$, $\sum_{k=1}^{\infty} k h_k < 1$, 
be a
subcritical offspring distribution for which there exists $c > 0$ such
that $\sum_{l \geq k }h_k \leq \exp \big\{ - ck \big\}$ for each $k \in \N$.  
We write $m_h = \sum_{k=1}^{\infty} k h_k$ for the mean number of offspring.
\end{hyp}

We now describe the law of the random biases.
\begin{definition}\label{deftreeweights}
Let $\nu$ be a probability distribution on $(1,\infty)$. For a tree T sampled from  $\P_h$, let $(\beta_e)_{e\in E(T)}$ be independent  and identically distributed random variables, with common distribution $\nu$. We write $ \mathbb{P}_{h,\nu}$ for the Galton-Watson tree with offspring distribution $h$ equipped with this set of random biases.
\end{definition}  

We make the following important non-lattice assumption on the distribution $\nu$.
\begin{hyp}\label{hypnu}
There exist $\bsup > \binf > 1$ such that the support of the measure $\nu$ is contained in $[\binf,\bsup]$.
Moreover, the support of $\nu \circ \log^{-1}$ is non-lattice. That is, the $\Z$-linear span of
$\log {\rm supp} (\nu)$ is dense in $\mathbb{R}$.
\end{hyp}
The fact that the distribution of the bias is supported in $(1,\infty)$ means that the walk is biased away from the root. We now need to introduce the following important exponent.
\begin{definition}
Let $\chi > 0$ be the unique value satisfying
\begin{equation}\label{formchi}
\int_1^{\infty} y^\chi    \nu (dy) = \frac{1}{m_h}.
\end{equation}
\end{definition}
Our first result shows that, under the assumptions above, the mean return time to the root has a pure power-law tail. The regularity of the tail of the random variable constrasts with the case of constant bias studied in \cite{BAFGH}.
\begin{definition}
Under the law $\P_{T,\beta}^\phi$, let $H_{\phi}= \inf \big\{ n \geq 1 , X_n=\phi \big\}$ denote the hitting time of the root.
\end{definition}
 Here and throughout, by $A(x) \sim B(x)$ is meant
$\frac{A(x)}{B(x)} \to 1$ as $x \to \infty$.
\begin{theorem}\label{theoremzero}
Assume Hypotheses \ref{hyph} and \ref{hypnu}. Then,
here exists a constant $\czerosm \in (0,\infty)$ such that
$$
 \mathbb{P}_{h,\nu} \Big(\E_{T,\beta}^{\phi}(H_{\phi})  > x \Big) \sim  \czerosm x^{-\chi}. 
$$ 
where  $\czerosm \in (0,\infty)$ is a constant that will be specified after Theorem \ref{theoremone}.
\end{theorem}

The electrical resistance theory of random walks, or, equivalently, the theory of random walk on weighted graphs (see \cite{Barlow} or \cite{Kumagai}), may be used to compute mean return times in terms of the total conductance (or weight) of the tree. Indeed, defining $\omega(T)$ to be this total weight $\sum_{e \in E(T)} c(e)$, we see that
\begin{equation}\label{defomegatree}
\omega(T) = \sum_{v \in V(T)} \omega(v). 
\end{equation}
The mean return time is then given by
\begin{equation}\label{ebetphi}
 \E_{T,\beta}^{\phi}(H_{\phi}) = \frac{2}{N} (\omega(T)-1)
\end{equation}
where $N$ be the number of offspring of the root. This formula arises by averaging the return time over the offspring of the root visited at the first step of the walk, with the commute-time formula of electrical resistance theory being used to express the mean return time from offspring to root. See the upcoming (\ref{eqelres}) for the relevant expression for this mean return time.

Using (\ref{ebetphi}), Theorem  \ref{theoremzero} is a simple consequence of the following pure power-law tail for the quantity $\omega(T)$.

\begin{theorem}\label{theoremone}
Assume Hypotheses \ref{hyph} and \ref{hypnu}.
There exists a constant $\cpo \in (0,\infty)$ such that
$$
 \mathbb{P}_{h,\nu} \Big( \omega(T) > u \Big) \sim \cpo u^{-\chi}, 
$$ 
\end{theorem}
 \noindent{\bf Remark.}
 The two constants $\czerosm$ and $\cpo$ in Theorems \ref{theoremzero} and \ref{theoremone} are explicitly given by
 \begin{equation}
 \czerosm= 2^{-\chi} \cpo \sum_{k=1}^{\infty}h_kk^{1-\chi}
 \end{equation}
 and, denoting by  $D(T) = \sup \big\{ d(\phi,v) : v \in V(T) \big\}$  the depth of the tree $T$,
\begin{equation}\label{formcpo}
 \cpo = 
\frac{1}{\chi m_h} \Big( \int_1^{\infty} y^\chi  \log(y)  \nu (dy) \Big)^{-1}
\lim_{k \to \infty} \E_{h,\nu} \Big( \omega(T)^{\chi} 1\!\!1_{D(T) = k}  \Big).
\end{equation}

It is natural to ask about not only the mean of the return time but also about how it is distributed.
We will present here an asymptotically accurate description of this distribution in the case that the tree is large (in a natural sense). This description plays a central role in the analysis of randomly biased walk on a supercritical Galton-Watson tree undertaken in \cite{GerardAlan}. We will see that the behaviour of the return time distribution has two entirely distinct regimes, according to whether or not the walker has the opportunity to visit the base of the tree (and so explore most of the tree) before its return to the root. We will condition on the events that the tree is large and that the walker makes such a deep visit, and conclude that the tail of the return time is then exponential, completing the analogy with the Bouchaud trap model that we have mentioned. 

Firstly, we introduce notation for conditioning that a tree be ``large'' enough. A tree $T$ will be regarded as large if $\omega(T)$ is large.
\begin{definition}\label{defphnu}
 For any $u > 0$ we set 
$\P_{h,\nu,u} = \P_{h,\nu} \big( \cdot \big\vert \omega(T) > u \big)$.
\end{definition}

We then introduce our notation for the distribution of the random walks. We sample a tree $T$ equipped with its biases $(\beta_e)$'s ,  choose a vertex $v \in V(T)$, and run the Markov chain on $V(T)$, started at $v$,
that is specified in Definition \ref{defwalk}.  
We denote by $\P_{h,\nu} \times \P_{T,\beta}^v$ the joint distribution of the tree, the biases and the random walk  started at $v$. We will also need to sample the weighted tree from the conditioned measure $\P_{h,\nu,u}$. Naturally,  we will then denote the joint distribution of the tree, the biases and the walk by $\P_{h,\nu,u} \times \P_{T,\beta}^v$.

A technical point is that, for a random tree $T$, we need to specify a selection rule for the starting point $v \in V(T)$ of the walk.
A selection rule can be easily defined using the classical  formal coding of  Galton-Watson trees (see  \cite{LeGall} Section 1.1).  In Appendix \ref{appa}, we review this coding and the definition of a selection rule.  

We also need to introduce notation for conditioning that the random walk explore enough of the tree.
In the next definition, we again use the coding of trees in Appendix \ref{appa} to specify a lexicographical ordering on the set of vertices of any Galton-Watson tree. 
\begin{definition}\label{vbase}
We call $\vbase$  the lexicographically minimal vertex of maximal depth in $V(T)$.
 \end{definition}
\begin{definition}
%Let us denote by $\mathcal{L}(T)$ the last generation of the tree $T$, i.e., the set of vertices of $T$ which are at maximal distance $D(T)$ from the root.
Under the law $\P_{T,\beta}^{\phi}$, the walk $X$ is said to make a {\it deep excursion} into $T$ if it hits $\vbase$ before returning to the root. We write $\falldeep = \big\{ H_{\vbase} < H_\phi \big\}$ for the event that $X$ makes a deep excursion into $T$. We will write $\P_{T,\beta, \falldeep}^{\phi}$ the distribution of the random walk conditioned by the event that it makes a deep excursion.
\end{definition}
We remark that, for a high choice of $u$, 
it is reasonable to suppose that, under the law
 $\P_{h,\nu,u}$, the tree is typically long and thin, consisting of a long path connecting the root to the base vertex $\vbase$, with only small subtrees hanging off this path. Indeed, in Theorem 
 \ref{lemdivexpold}, we will present a result to this effect.
 For a typical large tree, then, a walk from $\phi$ that makes a deep excursion into $T$ will be forced after reaching $\vbase$ to make a geometrically distributed number of highly ineffectual attempts to reach the root from $\vbase$, so that $H_\phi$ will be approximately exponentially distributed in this case. This is the content of our next result.
\begin{theorem}\label{thmtwoaah}
Assume Hypothesis \ref{hyph}, and that
there exist $\bsup > \binf > 1$ such that the support of the measure $\nu$ is contained in $[\binf,\bsup]$.
For all $t > 0$,
  $$
  \Big( \P_{h,\nu,u} \times \P_{T,\beta,\falldeep}^{\phi} \Big) \Bigg(  \frac{H_\phi}{\E_{T,\beta,\falldeep}^\phi(H_\phi)} > t \Bigg) \to \exp \big\{ -t \big\}
$$
as $u \to \infty$.
\end{theorem}
We also present an asymptotically accurate expression for the mean return time 
$\E_{T,\beta,\falldeep}^\phi(H_\phi)$. To do so, we introduce some notation.
\begin{definition}
Given a tree $T$ and $v \in V(T)$, we define the descendent tree $T_v$ of $v$ to be the subgraph induced by the set of all descendents of $v$ (among which, we include $v$). The root of $T_v$ is taken to be $v$.
\end{definition}
\begin{definition}\label{defomegachild}
Let $T$ be a weighted tree. Set $\vchild$ to be the neighbour of the root lying in the path $P_{\phi,\vbase}$ from the root to $\vbase$. The descendent tree $T_{\vchild}$ is itself a weighted tree with root $\vchild$.
We set $\omegachild$ equal to the weight of this tree, as defined by the formula (\ref{defomegatree}), where, to compute the summand on the right-hand-side, we note that, for each element $v \in V(T_{\vchild})$, the weight of $v$ is computed with respect to the root being $\vchild$. 
% $\sum_{v \in V(T_\vchild)} \omega_{\vchild}\big( v \big)$
%equal to the weight of the descendent tree $T_{\vchild}$ viewed from its root $\vchild$.

We further write $\pfd = \P_{T,\beta}^\vchild \big( \falldeep \big)$ for the probability that the walk starting at $\vchild$ makes a deep excursion into $T$.
\end{definition}
The significance of $\omegachild$ is its appearance in the mean hitting time formula
\begin{equation}\label{eqelres}
E_{T,\beta}^\vchild \big( H_\phi \big) = 2 \omegachild - 1,
\end{equation}
which follows from the commute-time formula in a reversible network, originally proved in \cite{CRRST}, and 
presented as Theorem 3.3 in \cite{Barlow}. 
\begin{theorem}\label{theoremthree}
Assume the hypotheses of Theorem \ref{thmtwoaah}.
For each $\epsilon > 0$,
$$
 \P_{h,\nu,u} \Bigg(   \frac{\E_{T,\beta,\falldeep}^\phi(H_\phi)}{ 2\omegachild \pfd^{-1}} \in \Big( 1 - \epsilon,1 + \epsilon \Big)  \Bigg) \to 1 
$$ 
as $u \to \infty$.

For all $t > 0$,
  $$
  \Big( \P_{h,\nu,u} \times \P_{T,\beta,\falldeep}^{\phi} \Big) \Bigg(  \frac{H_\phi}{2\omegachild \pfd^{-1}} > t \Bigg) \to \exp \big\{ -t \big\}
$$
as $u \to \infty$.
\end{theorem}
The new information presented in Theorem \ref{theoremthree} beyond that of Theorem \ref{thmtwoaah} is that,
 for a typical large tree, 
under $\P_{T,\beta}^\phi$, if $\falldeep$ does not occur, then $H_\phi$ is negligible. Indeed,
 that $H_\phi$ is negligible in the event $\falldeep^c$ is equivalent to saying that most of the mean $\E_{T,\beta}^\phi (H_\phi)$ arises on the event $\falldeep$. It easy to see that this in turn is the same as saying that 
$\E_{T,\beta,\falldeep}^\phi(H_\phi)$ is well approximated by 
 $2\omegachild \pfd^{-1}$. (The key ingredients to see the latter equivalence are (\ref{eqelres}), and the fact that, under $\falldeep$, the walk from $\phi$ must immediately travel to $\vchild$.)  
\begin{subsection}{The structure of the paper}
We begin in Section \ref{sectwo} by introducing some necessary tools from classical defective renewal theory and the extension needed for the proof of Theorem \ref{theoremone}, which forms the core of the argument for Theorem \ref{theoremzero}. We will illustrate how the defective renewal theory can be used straightforwardly to show that,  for a vertex $v \in V(T)$ of maximal depth, $\omega(v)$ has a pure power-law tail. This power-law differs by a constant from the one given in Theorem \ref{theoremone} for $\omega(T)$, the reason being that several vertices near the base of the tree contribute to $\omega(T)$ on roughly equal terms. In order to better approximate $\omega(T)$ and thus prove Theorem  \ref{theoremone}, we introduce a convenient decomposition of the tree in Section \ref{secthree}. The main result, Theorem \ref{lemdivexpold}, about this decomposition states that the components are small conditionally on $\omega(T)$ being large. The proof of Theorem \ref{lemdivexpold} is
fairly involved and in places delicate, and it is deferred to Appendix \ref{appb}. In Section \ref{secfour}, we prove Theorem \ref{theoremone}. In Section \ref{secfive}, we study the exponential tail of the return time, and prove Theorems \ref{thmtwoaah} and \ref{theoremthree}.

In its exploration of time spent in traps, the paper \cite{GerardAlan} requires a different decomposition from the one introduced in Section \ref{secthree}. In Appendix \ref{appc}, we introduce this new decomposition of the tree, which we call the renewal decomposition. The renewal decomposition has components that enjoy more independence than the corresponding objects in the decompositions used in this paper. 
The analogue  for the renewal decomposition of Theorem \ref{lemdivexpold} is needed in \cite{GerardAlan}.  Its proof also appears in Appendix \ref{appc}.
\end{subsection}
\end{section}

\begin{section}{A first use of defective renewal theory}\label{sectwo}

We recall and improve slightly on the classical non-lattice defective renewal theorem \cite{Feller}.  
Suppose that a succession of light bulbs has been produced. 
Independently of the earlier ones, a given light bulb is {\it permanent} with probability $1-p$, in which case, it will shine eternally. 
A light bulb which is not permanent is called {\it temporary}; given that a bulb is temporary (and independently of the status of the preceding bulbs), 
it has a lifetime distributed according to the law $\mu$.
At time zero, the first bulb is installed and begins to shine. A bulb is replaced whenever it fails, until a permanent one is installed.
The defective renewal theorem is concerned with the law of the time at which the final replacement is made 
(after which, a bulb will eternally shine). This is the time $Z$ in the next lemma, whose first part is the classical defective renewal theorem.
\begin{lemma}\label{lemabsfeller}
On a probability space $(\Omega,\P)$, let $\big\{ X_i: i \in \N  \big\}$ be a sequence of independent and identically distributed random variables.
The common distribution $\mu$ is assumed to have a non-negative, bounded and non-lattice support. 
Let $Y$ denote an independent
random variable taking values in the positive integers. 
 Set $Z = \sum_{i=1}^Y X_i$.
\begin{enumerate}
\item 
 Suppose that, under $\P$, 
$Y$ has the law of a geometric random variable of parameter $p \in (0,1)$: that is, for each
 $j \geq 0$, 
$\P \big( Y = j \big) = p^j (1-p)$.  Define $\kappa > 0$ as the unique value such that $\E(e^{\kappa X})= \frac{1}{p}$, and let $\conesm=\frac{1-p}{\kappa p \E(Xe^{\kappa X})}$.   
Then
$$
 \P \big( Z \geq u \big) \sim \conesm \exp \big\{ - \kappa u \big\},
$$
where recall that $f(u) \sim g(u)$ denotes $\lim_{u \to \infty} \frac{f(u)}{g(u)} = 1$. 
\item Suppose that there exist $p \in (0,1)$, 
$c \in (0,\infty)$ 
and 
$k_0 \in \N$ 
such that, for
$k \geq k_0$, 
$\P \big( Y = k \big) = c p^k$. Then 
$$
 \P \big( Z \geq u \big) \sim \ctwosm \exp \big\{ - \kappa u \big\},
$$
with $\ctwosm =  \frac{c\conesm}{1-p}$.
\item Suppose that, for such $p$ and $c$, we assume merely that
$\P \big( Y = k \big) \sim c p^k$ as $k \to \infty$. In this case also,  
$$
 \P \big( Z \geq u \big) \sim \ctwosm \exp \big\{ - \kappa u \big\}.
$$
\end{enumerate}
\end{lemma}
\noindent{\bf Proof.}
The first statement is (6.16) on page 377 of \cite{Feller}. (Note that there is typographical error in (6.16) and that $\mu$ should be replaced by $\mu^{\sharp}$.)

In the second statement, $Y$ assumes high values according to a law whose density is adjusted from that of a geometric law by multiplication by a factor of $\alpha = c/(1-p) \in (0,\infty)$.
Suppose that $\alpha \in (0,1]$. 
Then we may equip the probability space $\big( \Omega , \P \big)$
with
 a geometric random variable $V$ 
and an independent event $A$ satisfying $\P(A) = \alpha$ in such a way that
$Y 1\!\!1_{Y \geq k_0} = V 1\!\!1_{V \geq k_0} 1\!\!1_A$. 
Writing $Q$ for the supremum of the support of $\mu$, we have that,
if $u \geq k_0 Q$, then
 $Z \geq u$ implies that $Y \geq k_0$.
  We find then that, for such $u$, 
$\P( Y \geq u ) = \P ( V \geq u) \P(A) = \alpha \P (V \geq u)$, 
so that the second statement of the lemma follows from the first one in this case. 

In the case that $\alpha > 1$,    we similarly equip $\big( \Omega , \P \big)$ 
with a geometric random variable $V$ and an independent event
$A$ for which $\P(A) = \alpha^{-1}$ such that $V 1\!\!1_{V \geq k_0} = Y 1\!\!1_{Y \geq k_0} 1\!\!1_A$. 
We may then treat this case analogously to the preceding one.

The third statement may be reduced to the second one by a coupling argument. 
Let $Y'$ denote a random variable satisfying 
$\P \big( Y' = k \big) = c p^k$ for $k \geq k_0$, for some $k_0 \in \N$. 
We claim that a coupling $\coup$ of $Y$ and $Y'$ may be effected such that
\begin{equation}\label{ypreq}
\lim_{n \to \infty} \coup \Big( Y = Y' \Big\vert Y' = n \Big) = 1
\end{equation}
and
\begin{equation}\label{yeq}
\lim_{n \to \infty} \coup \Big( Y = Y' \Big\vert Y = n \Big) = 1.
\end{equation}
 Indeed, we may construct $\coup$ to satisfy $\coup \big( Y = n , Y' = n \big) = p_n \wedge q_n$,
with $p_n = \P ( Y = n)$ and $q_n = \P ( Y' = n )$. We then have that 
$\coup \big( Y \not= Y' \big\vert Y' = n \big) = \max \big\{ 0, 1 - p_n/q_n \big\}$,
 which tends to $0$ as $n \to \infty$. This yields (\ref{ypreq}), with the same argument giving (\ref{yeq}).

We now define $Z' = \sum_{i=1}^{Y'} X_i$, on $\big(\Omega,\P\big)$, 
and claim that
\begin{equation}\label{phbfracv}
 \lim_{u \to \infty} \frac{\mathbb{P} \big( Z' > u \big)}{\P  \big( Z > u \big)} = 1,
\end{equation}
so that the third statement of the lemma follows from the second applied to $Z'$. 
To prove (\ref{phbfracv}),
we construct under the measure $\coup$ the sequence $\big\{ X_i: i \in \N \big\}$ having the same law as under $\P$, and doing so independently of $Y$ and $Y'$. In this way, $\coup$ provides a coupling of $Z$ and $Z'$. 
It suffices for  (\ref{phbfracv})  to show that
\begin{equation}\label{bau}
\lim_{u \to \infty} \coup \Big( Z' > u \Big\vert Z > u \Big) = 1
\end{equation}
and
\begin{equation}\label{abu}
\lim_{u \to \infty} \coup \Big( Z > u \Big\vert Z' > u\Big)
   = 1.
\end{equation}
Given $Z \geq u$, it suffices for $Z' \geq u$ that $Y = Y'$. 
The conditional distribution of $Y$, given that $Z \geq u$, is supported on $[u/Q,\infty)$,
where recall that $Q = \sup {\rm supp}(\mu)$. Moreover, given $Y$, 
the event that $Y = Y'$ is conditionally independent of $Z$. Hence, 
$\coup \big( Y = Y' \big\vert Z \geq u \big) \to 1$ as $u \to \infty$ is a consequence of (\ref{yeq}).  
This proves (\ref{abu}). 
Likewise, (\ref{bau}) is derived by means of (\ref{ypreq}). $\Box$ \\

We now present Theorem \ref{theoremfive}, which states that $\omega(v)$  has a power-law tail for a vertex $v \in V(T)$ of maximal depth. The main tool is Lemma \ref{lemabsfeller}.
We will need the vertex $v$ to be  chosen independently of the biases $\{\beta_e: e \in E(T)\}$. For this theorem, any algorithm to do so would be satisfactory: for instance,  we could choose $v$ uniformly at random among those vertices of maximal depth, independently of the biases. For definiteness, we specify one such algorithm: indeed, we have already done so, in specifying the vertex $\vbase$ in Definition \ref{vbase} by means of  the classical  formal coding of the Galton-Watson tree (see  \cite{LeGall} Section 1.1, and Appendix \ref{appa}).
\begin{theorem}\label{theoremfive}
Assume Hypotheses \ref{hyph} and \ref{hypnu}.
The weight of the base vertex satisfies
 $$
 \mathbb{P}_{h,\nu} \Big( \omega(\vbase) > u \Big) \sim  \cpob u^{-\chi}, 
 $$
where $\chi$ is given in (\ref{formchi}) and where
\begin{equation}\label{formcpob}
 \cpob = \frac{\alpha}{\chi m_h \int_1^\infty y^{\chi} \log(y) d\nu(y)},
 \end{equation}
with the constant $\alpha \in (0,\infty)$ being given by the limit
\begin{equation}\label{alpha}
\alpha:= \lim_{n \to \infty} m_h^{-n} \P_h \Big( D(T) = n \Big)  \in  (0,\infty).
\end{equation}
\end{theorem}
The fact that  the distribution of the depth is asymptotically geometric, i.e., that the limit (\ref{alpha}) exists, is classical under  Hypothesis \ref{hyph}. See \cite{hsvj}. \\
\noindent{\bf Proof.}
We firstly note that the existence of the limit  (\ref{alpha}) is classical under  Hypothesis \ref{hyph}. See \cite{hsvj}.

The vertex $\vbase$ is at distance $D(T)$ from the root $\phi$.
Under $\P_{h,\nu}$, the biases attached to the path from
$\phi$ to $\vbase$ are independent samples of the law $\nu$.
Hence, under $\P_{h,\nu}$,
\begin{equation}\label{logomid}
\log \omega(\vbase) =  \sum_{i=1}^{D(T)} X_i,
\end{equation}
where $\big\{ X_i: 1 \leq i \leq D(T) \big\}$ is a sequence of independent random
variables, each having the law $\nu \circ \log^{-1}$.
Note that this sum satisfies the hypotheses of the third part of Lemma \ref{lemabsfeller}, due to $D(T)$ being independent of $\big\{ X_i: i \in D(T) \big\}$, and (\ref{alpha}). 
Applying this result, we find that
\begin{equation}\label{simlogom}
\mathbb{P}_{h,\nu} \Big( \log \omega(\vbase) > u \Big) \sim \cpob e^{- \chi u}.
\end{equation}
where $\chi > 0$ satisifes (\ref{formchi})
and where $\cpob$ has the form given in (\ref{formcpob}). $\Box$ 

The quantity $\omega(\vbase)$ is a first approximation to $\omega(T)$, for a tree sampled under the measure $\P_{h,\nu,u}$. This approximation is never accurate up to leading order, although it is plausible that, for a typical sample of this measure,  $\omega(T)$ is approximated up to a small multiplicative correction by a sum of $\omega(v)$ over vertices $v$ ranging over several generations near the end of the trap. 
We will introduce such an approximation, and analyze it, also by means of the defective renewal theorem. We now define a decomposition of a tree, in order to make this approximation. 
 \end{section}
 
\begin{section}{Bare trees and approximation of the total weight}\label{secthree}
\begin{subsection}{Approximating $\omega(T)$ using a splitting of the tree}

 We consider firstly a simple splitting of the tree $T$ and a convenient approximation of 
 its total weight $\omega(T)$. 
 \begin{definition}\label{defoutgr}
$\empty$
 \begin{enumerate}
 \item We record the vertices in $P_{\phi,\vbase}$ in the form $\big[ \phi = \psi_0,\psi_1,\ldots,\psi_{D(T)} = \vbase \big]$.
\item For $0 \leq i \leq D(T)$, let $J_i$ denote the connected component containing $\psi_i$ of the graph with vertex set $V(T)$ and edge-set 
$E(T) \setminus E \big( P_{\phi,\vbase} \big)$. We will call  $J_i$  the $i$-th outgrowth of $T$. Note that $J_{D(T)}$ is the singleton graph with vertex $\vbase$.
\item For any integer $k$, define $w_k(T) = \sum_{i=0}^k \sum_{v \in V(J_{D(T)-i})} \omega(v)$.
\end{enumerate}
 \end{definition}
The quantity $w_k(T)$
will provide a better approximation to $\omega(T)$ than that which $\omega(\vbase)$ offers. Before explaining how this is so, we need a little more notation.
%\begin{definition}
%Given a tree $T$ and $v \in V(T)$, we define the descendent tree $T_v$ of $v$ to be the subgraph induced by the set of all descendents of $v$. The root of $T_v$ is taken to be $v$.
%\end{definition}
\begin{definition}\label{defek}
 Let $k \in \N$. Let $T$ be a weighted tree for which $D(T) \geq k$.
Set $E_k$ to be the weighted descendent tree $T_{\psi_{D(T)-k}}$, or, equivalently, 
the weighted subgraph induced by the set of 
 vertices $v \in \cup_{i = D(T) - k}^{D(T)} V(J_i)$. 
%lying in one of the $(k+ 1)$-st final outgrowths of $T$. 
Note that $w_k(T) = \sum_{v \in V(E_k)} \omega(v)$. 

We also set $E_k^*$ equal to the weighted subtree of $T$ induced by  the vertex set given by removing the strict descendents of
$\psi_{D(T)-k}$ from $V(T)$. Note that $E\big( E_k \big) \cup  E\big( E^*_k \big)$ is a partition of $E(T)$. 
\end{definition}
We extend the notation $\omega(v)$ for the weight of a vertex $v \in V(T)$
in a tree $T$ in the following way.
\begin{definition}\label{defomegaext}
Let T be a weighted tree, and let $u,v \in V(T)$, with $v$  being a descendent of $u$.
We write $\omega_u(v)$ for the product of the edge-weights in the path $P_{u,v}$. Note that, for any $u \in V(T)$,
we have that $\omega_\phi(u) = \omega(u)$.
\end{definition}
The quantity $w_k(T)$ has a convenient representation. 
For any $v \in V(E_k)$, we may write $P_{\phi,v}$ as the concatenation of $P_{\phi,\psi_{D(T) - k}}$
and  $P_{\psi_{D(T) - k},v}$. By so doing, we obtain
$\omega (v) = \omega \big( \psi_{D(T) - k} \big) \omega_{\psi_{D(T) - k}}(v)$. 
By summing this formula over $v \in V(E_k)$, we arrive at
\begin{equation}\label{eqprod}
w_k(T) = u_k(T) v_k(T).
\end{equation}
Here, $u_k(T) : = \omega \big( \psi_{D(T) - k} \big)$ and 
$v_k(T) : = \sum_{v \in V(E_k)} \omega_{\psi_{D(T) - k}}(v)$. This product formula will be essential for our analysis of the approximation $w_k(T)$.

\end{subsection}
\begin{subsection}{Outgrowths in a high-weight tree are small}
To obtain Theorem \ref{theoremone}, as well as deriving the asymptotic behaviour of 
$w_k(T)$, we must show that $w_k(T)$ is a good enough approximation to $\omega(T)$.
 To do so, we need to know that the,
outgrowths in a typical sample of $\P_{h,\nu}$ are not too large, 
uniformly under conditioning on the sample being a tree of any high weight. 
This estimate is provided by the next theorem, which is stated using the notation introduced in Definition \ref{defphnu}.
\begin{theorem}\label{lemdivexpold}
For the statement, we take $J_i = \emptyset$ if $i > D(T)$ (for any 
tree $T$). 
Let $\big\{ h_i: i \in \N \big\}$ satisfy Hypothesis \ref{hyph}, and let $\nu$ be an edge-weight law with compact support in $(1,\infty)$.
There exists $\conlem > 0$ such that, for all $u > 0$ and $i \in \mathbb{N}$,
$$
\mathbb{P}_{h,\nu,u} \Big( \big\vert V \big( J_i \big) \big\vert \geq k \Big) \leq \exp \big\{ - \conlem k  \big\},
$$
for each $k \in \mathbb{N}$.
\end{theorem}
While the statement of Theorem \ref{lemdivexpold} is natural enough, its proof is a little technical. We defer it to Appendix \ref{appb}.

We also record the following property of the decomposition.
\begin{lemma}\label{lemekprop}
The law of the tree $E_k$ under $\P_{h,\nu} \big( \cdot \big\vert D(T) \geq k \big)$
coincides with that of $T$ under $\P_{h,\nu} \big( \cdot \big\vert D(T) = k \big)$. Moreover, this statement holds if 
$\P_{h,\nu} \big( \cdot \big\vert D(T) \geq k \big)$ is conditioned on any 
admissible choice of $E^*_k$.
\end{lemma}
\noindent{\bf Proof.} 
Suppose that a sample $T$ of  
$\P_{h,\nu} \big( \cdot \big\vert D(T) \geq k \big)$ 
 is further conditioned on an admissible choice of $E_k^*$. (We write $\tilde\P$ for the conditioned law.)
The tree $E_k$, which, by definition, is given by the descendent tree $T_{\psi_{D(T) - k}}$, must have $D (E_k) = k$, by the definition of $\psi_{D(T) - k}$. Moreover, $E_k$ may take the value of any such tree. To see this, 
note firstly that 
there must exist a tree $T'$ with $D(T') = k$  
in the support of the random variable $E_k$  under $\tilde\P$, 
since the choice of $E_k^*$ under which we condition to obtain $\tilde{\P}$ is 
an admissible one; 
and note further that the tree obtained by making the choice $E_k = T'$  has $\vbase \in V(E_k)$. 
This forces $\psi_{D(T) - k}$ to be lexicographically smaller than any vertex $v \in V(T) \setminus V(E_k)$  for which $d\big(\phi, v   \big) = D(T)$. This means that, for an arbitrary choice of $E_k = T''$, with $T''$ a tree for which $D(T'') = k$, we necessarily have $\vbase \in V(E_k)$, so that, indeed, $E_k$ may take the value $T''$. 

Hence, under $\tilde\P$, $E_k$ has the law of $\P_{h,\nu}$ subject merely to the condition that $D (E_k) = k$.
$\Box$
\end{subsection}

\begin{subsection}{Bare trees}
The following definition and lemma reformulate Theorem \ref{lemdivexpold} in a manner
convenient to applications.
\begin{definition}\label{defbt}
For $\barecon > 0$, we say that a  \wgt tree $T$ is $\barecon$-bare if
$\big\vert V(J_i) \big\vert \leq \barecon \log \log \omega(T)$
for each $i \in \{ 0, \ldots , D(T) - 1 \}$.
\end{definition}
\begin{lemma}\label{lemcreg}
Let $\conlem > 0$ denote the constant that appears in Theorem \ref{lemdivexpold}.
Fix $\barecon > 2/\conlem$. Then
$$
\mathbb{P}_{h,\nu,u} \Big( T \, \textrm{is $\barecon$-bare}   \Big) \geq 1 - (\log u)^{- \conlem \barecon/2},
$$
for sufficiently high $u$.
\end{lemma}
\noindent{\bf Proof.} 
Note that $\omega(T) \geq \omega(\vbase) \geq \binf^{D(T)}$. Hence, by  (\ref{alpha}),

\begin{equation}\label{omvbdold}
  \P_{h,\nu} \Big( \omega(T) > u \Big) \geq \P_{h,\nu} \Big( D(T) > \frac{\log u}{\log \binf} \Big)
  \sim \cone u^{- \frac{\log\big(m_h^{-1}\big)}{\log \binf}},
\end{equation}
for some $\cone \in (0,\infty)$.

By (\ref{alpha}) once more, and (\ref{omvbdold}), we see that, for any $C > \frac{1}{\log \binf}$ and for all sufficiently high $u$, 
\begin{eqnarray}
 & & \mathbb{P}_{h,\nu,u}\Big(  D(T) > C \log u \Big)
 \leq   \frac{\mathbb{P}_h \big( \vert D(T) \vert \geq C \log u
   \big)}{\mathbb{P}_{h,\nu} \big( \omega(T) > u \big) } \nonumber \\
 & \leq &  
 2 u^{- \Big( 1/{\big( \log \binf \big)}  - C  \Big) \log\big(m_h^{-1}\big)}. \nonumber
\end{eqnarray}

We find then that
\begin{eqnarray}
 & & \mathbb{P}_{h,\nu,u} \Big( T \, \textrm{is not $\barecon$-bare}   \Big) 
                               \nonumber   \\
 & \leq & \mathbb{P}_{h,\nu,u} \Big(  \max_{i \in \{ 0,\ldots,D(T) - 1\} } \vert V(J_i) \vert > \barecon \log \log u \Big) \nonumber \\
 & \leq & C \big( \log u \big) \exp \big\{ - \conlem \barecon \log \log u \big\} \, + \,
  \mathbb{P}_{h,\nu,u} \big(  D(T) > C \log u \big) \nonumber \\
  & \leq &  C \big( \log u \big)^{1- \conlem \barecon} \, + \, 
 2 u^{- \Big( 1/{\big( \log \binf \big)}  - C  \Big) \log\big(m_h^{-1}\big)}, \nonumber
\end{eqnarray}
the second inequality by virtue of Theorem \ref{lemdivexpold}. 
Note that $\conlem \barecon > 2$ implies that $1 - \conlem \barecon \leq - (\conlem/2) \barecon$. This
yields the result. $\Box$ \\
We also record some properties of a $\barecon$-bare tree.
\begin{lemma}\label{lemom}
Let $T$ be a $\conbare$-bare tree.
There exists a constant $\conom > 0$ such that
$$
\omega\big( \vbase \big) \geq \frac{\omega(T)}{\big( \log \omega(T )\big)^{2 \conbare \log \bsup}}.
$$
and
$$
D(T) \geq \frac{\log \omega(T)}{2 \log \bsup},
$$
provided that $\omega(T) \geq \conom$.
\end{lemma}
\noindent{\bf Proof.}
Let $v \in V(J_k)$ for some $k \in \big\{ 0,\ldots,D(T)-1 \big\}$. 
Note that $\psi_k \in V(T)$ is the latest common ancestor of $v$ and $\vbase$.
Recalling that $d(\cdot,\cdot)$ denotes  the graphical distance on $V(T)$, note that
$$
\omega(v) \leq \bsup^{d\big(\psi_k,v\big)} \omega(\psi_k) 
\leq \bsup^{\conbare \log\log \omega(T)} \omega(\psi_k),
$$
where we used that $T$ is $\conbare$-bare in the second inequality, while
$$
 \omega(\vbase) \geq \binf^{d\big(\psi_k,\vbase\big)} \omega(\psi_k) = \binf^{D(T) - k} \omega(\psi_k).
$$
Thus,
$$
\omega(v) \leq \bsup^{\conbare \log \log \omega(T)} \binf^{-\big( D(T) - k \big)}
 \omega(\vbase).
$$
We find that
\begin{eqnarray}
 & & \omega(T) = \sum_{v \in V(T)} \omega(v)
 = \sum_{j=0}^{D(T)} \sum_{v \in V \big( J_{D(T)  - j} \big)} \omega(v) \nonumber \\
 & \leq & \omega(\vbase) \bsup^{\conbare \log \log \omega(T)}
 \sum_{k=0}^{D(T)} \binf^{-k} \big\vert V \big( J_{D(T) - k} \big) \big\vert \leq 
 \omega(\vbase) \bsup^{2\conbare \log \log \omega(T)}, \nonumber
\end{eqnarray}
the second inequality by $T$ being $\conbare$-bare and $\omega(T) > \conom$.
Thus, the first statement of the lemma.
Using $\omega(\vbase) \leq \bsup^{D(T)}$, we obtain the second. $\Box$
\end{subsection}
\end{section}

\begin{section}{The proof of Theorem \ref{theoremone}}\label{secfour}
To prove Theorem \ref{theoremone}, we will use the approximation of the weight $\omega(T)$ provided by $w_k(T)$ that was introduced in Definition \ref{defoutgr}. 
We will analyze this quantity by means of the product form (\ref{eqprod}).
We will prove that, under the law $\P_{h,\nu} \big( \cdot \big\vert D(T) \geq k \big)$, the two constituents $u_k(T)$ and $v_k(T)$ are
independent, and $u_k(T)$ is typically the larger.  As we will shortly explain, the asymptotic decay of $u_k(T)$ 
may be obtained by the defective renewal theory with which we obtained that of $\omega(\vbase)$ under $\P_{h,\nu}$ in the proof of Theorem \ref{theoremfive}.
We will see that the term
$v_k(T)$ modifies the determined decay rate only by a
constant. Indeed, 
this is reflected in the statement of
Theorem \ref{theoremone}, where the constant $\cpo$ is expressed as a limit of a certain average of vertex weights near the end of the tree.  

We begin by computing the tail of the main term $u_k(T)$ by another use of defective renewal theory. 
\begin{lemma}\label{lemchipr}
There exists a sequence of finite and positive constants $\big\{ \constd(k): k \in \N \big\}$ such that, for each $k \in \N$,
$$
\mathbb{P}_{h,\nu} \Big( u_k(T) \geq u \Big\vert D(T) \geq k \Big) \sim
 \constd(k) u^{-\chi}.
$$
We have that

$$
\lim_{k \to \infty} \constd(k) =   \frac{1 - m_h}{\chi m_h  \int_1^{\infty} y^{\chi} \log(y) d\nu(y)}.
$$ 
\end{lemma}
\noindent{\bf Proof.} 
Analogously to (\ref{logomid}), we have that, under $\P_{h,\nu} \big( \cdot \big\vert D(T) \geq k \big)$,
$$
 \log u_k(T) = \sum_{i=1}^{D(T) - k} X_i,
$$
where $\big\{ X_i: i \in \N \big\}$ is a sequence of independent random
variables, each having the law $\nu \circ \log^{-1}$.
The depth $D(T)$ is independent of $\big\{ X_i: i \in \N \big\}$ and, for 
every $k \in \N$, the following limit exists by (\ref{alpha}):
\begin{equation}\label{eqalk}
\lim_{n \to \infty}
m_h^{- n} \P_h \Big( D(T) = n + k \Big\vert D(T) \geq k \Big)  =: \alpha_k.
\end{equation}
as $n \to \infty$.
It follows directly from (\ref{eqalk}) that 
\begin{equation}\label{alpklim}
\lim_{k \to \infty} \alpha_k = 1 - m_h. 
\end{equation}
By the third part of Lemma \ref{lemabsfeller}, (and analogously to (\ref{simlogom})),
$$
  \P_{h,\nu} \Big(  \log u_k(T)  > u \Big\vert D(T) \geq k \Big) \sim 
   \constd(k) e^{-\chi u},
$$
as $u \to \infty$,
with
 $\constd(k) = \frac{\alpha_k}{\chi m_h \int_1^\infty y^\chi \log(y) d\nu(y)}$.

The stated convergence of $\big\{ \constd(k): k \in \N \big\}$   follows from (\ref{alpklim}). $\Box$ \\
We will now show that the two terms in the product decomposition (\ref{eqprod}) 
for $w_k(T)$
are independent, and that $u_k(T)$ is the dominant one:
\begin{lemma}\label{lrgone} Under the measure 
$\mathbb{P}_{h,\nu} \big( \cdot  \big\vert D(T) \geq k \big)$,
the random variables $u_k(T)$ and $v_k(T)$ are independent.
\end{lemma}
\noindent{\bf Proof.}
On the space 
$\mathbb{P}_{h,\nu} \big( \cdot  \big\vert D(T) \geq k \big)$,
the random variables $u_k(T)$ and $v_k(T)$ are respectively
measurable with respect to the labelling with weights of the edge-disjoint trees
$E^*_k$ and 
$E_k$ that were introduced in Definition \ref{defek}. The statement thus follows from Lemma \ref{lemekprop}. $\Box$
\begin{lemma}
There exist $c >0$ and $u_0:\N \to (0,\infty)$ such that, for any $k \in \mathbb{N}$, and for $u > u_0(k)$,
\begin{equation}\label{uepst}
 \mathbb{P}_{h,\nu} \Big(  v_k(T) > u \Big\vert D(T) \geq k \Big)
 \leq  \exp \big\{ - c u \big\}
 \mathbb{P}_{h,\nu} \Big(  u_k(T) > u \Big\vert D(T) \geq k \Big).
\end{equation}
\end{lemma}
\noindent{\bf Proof.}
Note that
\begin{equation}\label{eqvkt}
v_k(T) \leq \big\vert V(E_k) \big\vert \bsup^k.
\end{equation}
Note further that
\begin{eqnarray}
 & &   \P_{h,\nu} \Big(  \big\vert V(E_k) \big\vert \geq v \Big\vert D(T) \geq k \Big)
   =   \P_{h,\nu} \Big(  \big\vert V(T) \big\vert \geq v \Big\vert D(T) = k \Big) \nonumber \\
    & \leq & \frac{\P_{h,\nu} \Big(  \big\vert V(T) \big\vert \geq v  \Big)}{\P_{h,\nu} \Big(  D(T) = k  \Big)}
     \leq C_k \exp \big\{ - c v \big\} \label{eqvek}
 \end{eqnarray} 
Here, the equality is due to 
Lemma \ref{lemekprop}, and the second inequality, to the following Lemma \ref{sizedec}.

By (\ref{eqvkt}) and (\ref{eqvek}), 
$\mathbb{P}_{h,\nu} \big( v_k(T)  \geq u \big\vert D(T) \geq k \big) \leq \exp \big\{ - c\bsup^{-k} u \big\}$ for sufficiently high $u$. By Lemma \ref{lemchipr}, we obtain the statement.  $\Box$ 
\begin{lemma}\label{sizedec}
There exists a constant $c \in (0,1)$ such that
$$
\mathbb{P}_h \Big( \big\vert V(T) \big\vert \geq n \Big) \leq c^n 
$$
for each $n \in \mathbb{N}$.
\end{lemma}
\noindent{\bf Proof.} It is a classical fact (see \cite{LeGall}, Corollary 1.6) that the random variable $\vert V(T) \vert$ has the same distribution as the hitting time $ U= \inf \big\{ n \leq 1, S_n= -1 \big\}$, 
where $S_n= \sum_{i=1}^n (M_i-1)$  and the $\{ M_i: i \in \N \}$ are independent and identically distributed random variables with distribution $h$. The distribution of the increments $M_i$ has a negative mean and an exponential tail; thus, by a classical argument, the hitting time $U$ also has an exponential tail. $ \Box$

We are now ready to prove a precise decay rate for $w_k(T)$ under $\P_{h,\nu}$.
\begin{lemma}\label{lrgthr}
For each $k \in \mathbb{N}$, 
\begin{eqnarray}
 & & \mathbb{P}_{h,\nu} \Big(  \Big\{ w_k(T) > u \Big\} \cap \Big\{
 D(T) \geq k \Big\} \Big) \nonumber \\
 & \sim & \constd(k) \E_{h,\nu} \Big(   \omega(T)^{\chi} \Big\vert D(T) = k  \Big)  
   \P_{h,\nu} \Big( D(T) \geq k \Big)  u^{- \chi}, \nonumber 
\end{eqnarray}
where the constants $\big\{ \constd(k): k \in \N \big\}$ appear in the statement of Lemma
\ref{lemchipr}.
\end{lemma}
\noindent{\bf Proof.}
It suffices to show that, for each $k \in \mathbb{N}$,
\begin{equation}\label{omkasymp}
\mathbb{P}_{h,\nu} \Big( w_k(T) > u \Big\vert D(T) \geq k \Big)
 \sim  \constd(k) \E_{h,\nu} \Big( \big( v_k(T) \big)^{\chi} \Big\vert D(T) \geq k \Big) u^{-\chi}, 
\end{equation}
since Lemma \ref{lemekprop} demonstrates that the expectation on the right-hand-side is equal to 
$\E_{h,\nu} \big( \omega(T)^\chi \big\vert D(T) = k \big)$.

By Lemmas \ref{lemchipr} and \ref{lrgone}, (\ref{omkasymp}) follows from the next lemma. 
$\Box$
\begin{lemma}
Let $U$ and $V$ be independent random variables on a probability space
$(\Omega,\mathbb{P})$ such that
\begin{equation}\label{uuchi}
 \mathbb{P}\big(U > u\big) \sim c u^{-\kappa}
\end{equation}
for some $c > 0$ and $\kappa > 0$, and, for some $\eta > 0$, 
\begin{equation}\label{vuchi}
 \mathbb{P}\big( V > u \big) \leq u^{-\eta}
 \mathbb{P}\big( U > u \big) 
\end{equation}
for $u$ sufficiently high. Assume also that $V \geq 1$.
Then
\begin{equation}\label{uvchi}
 \mathbb{P}\big( U V > u \big) \sim c \E \big( V^{\kappa} \big) u^{-\kappa}. 
\end{equation}
\end{lemma}

\noindent{\bf Proof.}
 
Fix $\epsilon > 0$.
We choose $u_0 = u_0 (\epsilon)$
such that, for $u \geq u_0$,
\begin{equation}\label{uraeps}
 \Big\vert \frac{\P \big( U > u \big)}{c u^{-\kappa}}  - 1  \Big\vert < \epsilon.
\end{equation}
From (\ref{uuchi}), (\ref{vuchi}) and $V \geq 1$, it follows that, for $\delta > 0$,
there exists $C = C(\delta) > 0$ with $C(\delta) \to \infty$ as $\delta \to
0$ and $u_1 = u_1(\delta) > 0$ such that $u > u_1$ implies that
\begin{equation}\label{cvbn}
\mathbb{P} \big( UV > u \big) \leq \big( 1 + \delta \big)
 \mathbb{P} \Big( \Big\{ UV > u \Big\} \cap \Big\{ V \leq u/C \Big\} \Big).
\end{equation}
Now, choose $\delta > 0$ such that
$C(\delta) > u_0(\epsilon)$ and $\delta < \epsilon$.

Let $g$
denote the distribution function of $V$.
In what follows, each quantity $E_i$ is an error term, in absolute value at
most $\epsilon_1$: for $u \geq \max \big\{ u_0,u_1 \big\}$,
\begin{eqnarray}
 & &  \mathbb{P} \Big( UV > u  \Big) \nonumber \\ 
 & = & \big( 1 + E_1 \big) \mathbb{P} \Big( \Big\{ UV > u \Big\} 
  \cap \Big\{  V \leq  \frac{u}{C(\delta)} \Big\} \Big) \nonumber \\
 & = & \big( 1 + E_1 \big) \int_0^{\frac{u}{C(\delta)}}
 \mathbb{P} \Big(  U > \frac{u}{v} \Big) dg (v) \nonumber \\
 & = &  \big( 1 + E_1 \big) \big( 1 + E_2 \big) c u^{- \kappa}
\int_0^{\frac{u}{C(\delta)}} v^{\kappa} dg(v) \nonumber 
\end{eqnarray}
where, in the third equality, 
we used $u/v > u_0$ for $v \leq u/C(\delta)$ (which is implied by $C(\delta) >
u_0$), (\ref{uraeps}), and $\delta < \epsilon$. Noting that
$\lim_{u \to \infty}\int_0^{\frac{u}{C(\delta)}} v^{\kappa} dg(v)  = \E \big(
V^{\kappa} \big)$, and  $\E \big(
V^{\kappa} \big) < \infty$ by (\ref{uuchi}) and (\ref{vuchi}),
we obtain (\ref{uvchi}). $\Box$ \\

The next step is to show that the approximation of $\omega(T)$ by $w_k(T)$ is sufficiently good. 
\begin{lemma}\label{lrgtwo}
There exists $\epsilon(k) \in (0,\infty)$ with $\epsilon(k) \to 0$ as $k \to \infty$ and $u_0(k) \in (0,\infty)$ such that, for $u \geq u_0(k)$,
\begin{eqnarray}
 & & 
 \mathbb{P}_{h,\nu} \Big( \Big\{  w_k(T) > u  \Big\} \cap \Big\{ D(T) \geq k \Big\} \Big) \leq \mathbb{P}_{h,\nu} \Big( \omega(T) > u \Big) \nonumber \\
 & \leq & \Big( 1 + \epsilon (k) \Big) 
 \mathbb{P}_{h,\nu} \Big(  \Big\{ w_k(T) > u \Big\} \cap \Big\{ D(T) \geq k \Big\} \Big). \nonumber
\end{eqnarray}
\end{lemma}
\noindent{\bf Proof.}
The first inequality is implied by $w_k(T) \leq \omega(T)$. 
Regarding the second inequality, note that there exists $c > 0$ such that, for any $k \in \N$, and for $u > u_0(k)$,
\begin{equation}\label{dku}
\P_{h,\nu,u} \Big( D(T) \leq k \Big) \leq \exp \big\{ - c \bsup^{-k}u \big\}.
\end{equation}
Indeed, since $\omega(T) \leq \vert V(T) \vert \bsup^{D(T)}$, 
on $\big\{ \omega(T) > u, D(T) \leq k \big\}$, we have that $\vert V(T) \vert \geq u \bsup^{-k}$. By means of Lemma \ref{sizedec}, we have then that $$\P_{h,\nu} \big(  \omega(T) > u, D(T) \leq k \big) \leq \exp\big\{ - c u \bsup^{-k} \big\}$$ However, by $\omega(T) \geq \omega \big( \vbase \big)$ and Theorem \ref{theoremfive}, we have that $\P_{h,\nu} \big(  \omega(T) > u \big) \geq (\cpob/2 )u^{-\chi}$. Hence, we have (\ref{dku}).

In view of 
Lemma \ref{lrgthr} and (\ref{dku}), it is enough for the second inequality in the statement to show that,
 for all $\epsilon > 0$, there exists a $k_0(\epsilon)$ such that for $ k > k_0 $,
\begin{equation}\label{eqphueps}
\lim_{u \to \infty} \P_{h,\nu,u} \Big(  w_k(T) > u(1-\epsilon) \Big) =1.
\end{equation}
It is sufficient for (\ref{eqphueps}) that
\begin{equation}\label{phomomk}
\lim_{u \to \infty} \P_{h,\nu,u} \Big(  \omega(T) - w_k(T) \leq  \epsilon \omega(T) \Big) = 1.
\end{equation}
To show this, we begin by noting that $\omega(T) - w_k(T) = \sum_{i=0}^{D(T) - k - 1} \sum_{v \in V(J_i)} \omega(v)$.

Recall Definition \ref{defomegaext} and so note also that 
$$
\omega(T) \geq \omega(\vbase) = \omega \big( \psi_{D(T) - k} \big) 
\omega_{\psi_{D(T) - k}}(\vbase) \geq u_k(T) \binf^k.
$$ 
Hence, 
\begin{equation}\label{omineq}
  \frac{\omega(T) - w_k(T)}{\omega(T)}
\leq
 \binf^{-k} u_k(T)^{-1}  \sum_{i=0}^{D(T) - k - 1} \sum_{v \in V(J_i)} \omega(v).
\end{equation}
Let $j \in \{ 0 ,\ldots,D(T)-k - 1 \}$ and $v \in V \big( J_{D(T) - k - 1 - j} \big)$. 
We claim that
\begin{equation}\label{omuk}
  \omega(v) \leq u_k(T) \bsup^{\big\vert V \big( J_{D(T) - k - 1 - j} \big) \big\vert} \binf^{- j-1}.
\end{equation}
To prove this, note that
$\omega(v) = \omega\big( \psi_{D(T) - k - 1 - j} \big) \omega_{\psi_{D(T) - k - 1 - j}}(v)$,
while 
$$
u_k(T) = \omega \big( \psi_{D(T) - k} \big) = 
 \omega\big( \psi_{D(T) - k - 1 - j} \big) \omega_{\psi_{D(T) - k - 1 - j}}\big( \psi_{D(T) - k} \big).
$$
Hence,
\begin{equation}\label{omukquot}
  \frac{\omega(v)}{u_k(T)}   = \frac{\omega_{\psi_{D(T) - k - 1 - j}}(v)}{\omega_{\psi_{D(T) - k - 1 - j}}\big( \psi_{D(T) - k} \big)}. 
\end{equation}
Now, $\omega_{\psi_{D(T) - k - 1 - j}}(v) \leq \bsup^{\big\vert V \big( J_{D(T) - k - 1 - j} \big) \big\vert}$,
and
$$
\omega_{\psi_{D(T) - k - 1 - j}}\big( \psi_{D(T) - k} \big) \geq \binf^{j+1},
$$
so that (\ref{omuk}) follows from (\ref{omukquot}).

From (\ref{omuk}), we find that
$$
 \frac{1}{u_k(T)} \sum_{i = 0}^{D(T) - k - 1} \sum_{v \in V(J_i)} \omega(v)
  \leq 
 \sum_{i=0}^{D(T) - k - 1} \binf^{-1-i} \big( 2 \bsup \big)^{\big\vert
   V \big( J_{D(T) - k - 1 - i} \big) \big\vert}. \nonumber
$$
In seeking to verify (\ref{phomomk}),
we require a bound on the upper tail under $\P_{h,\nu,u}$ of the right-hand-side of the preceding inequality, which amounts, in effect, to an assertion that a tree sampled under $\P_{h,\nu,u}$ is typically sparse at a short distance from its end. 
The next lemma provides such a bound. It does not follow from the statement of Theorem \ref{lemdivexpold}, since the indices $D(T),D(T)-1,\ldots$ of the final outgrowths of $T$ are random. 
However, it is, in essence, a byproduct of the method of proof of Theorem \ref{lemdivexpold}, and its proof appears at the end of Appendix \ref{appb}.
\begin{lemma}\label{lemaux}
There exists a constant $c > 0$ such that, for all $\ell \in \N$,
$$
 \sup_{u > 0} \P_{h,\nu,u} \Big( \sum_{i=0}^{D(T) - \ell} \big(2\bsup\big)^{\big\vert V \big( J_{D(T) - \ell - i} \big) \big\vert}  
  \binf^{-1-i-\ell}  \geq \exp \big\{ - c\ell \big\} \Big) \leq \exp \big\{ - c\ell \big\}.
$$
\end{lemma}
Applying Lemma \ref{lemaux} with $\ell = k + 1$ in the second inequality below, we obtain
\begin{eqnarray}
 & & \P_{h,\nu,u} \Big( \frac{1}{u_k(T)} \sum_{i = 0}^{D(T) - k - 1} \sum_{v \in V(J_i)} \omega(v) > \exp\{ - c k \} \binf^{k} , D(T) \geq k \Big) \label{phuvineq} \\
 & \leq &
 \P_{h,\nu,u} \Big(
 \sum_{i=0}^{D(T) - k - 1} \binf^{-i - k} \big( 2 \bsup \big)^{\big\vert
   V \big( J_{D(T) - k - 1 - i} \big) \big\vert} > \exp \big\{ - c k \big\}
 \Big) \leq \exp \big\{ - ck \big\}.
 \nonumber
\end{eqnarray} 
In light of (\ref{omineq}) and (\ref{phuvineq}),
$$
\P_{h,\nu,u} \Big(  \omega(T)  - w_k(T) > \exp \{ - ck \}  \omega(T),  D(T) \geq k \Big)
 \leq \exp \{ - ck \}.
$$
Using (\ref{dku}), and letting $u \to \infty$ and taking $k$ large enough, we obtain (\ref{phomomk}). This completes the proof. $\Box$ \\
The following result, whose proof is trivial, is required to obtain 
Theorem \ref{theoremone}
from Lemmas \ref{lrgthr} and \ref{lrgtwo}.
\begin{lemma}\label{lrgfour}
Let $\gamma > 0$.
Let $s:[0,\infty) \to [0,\infty)$, 
$\big\{ s_{\epsilon}: [0,\infty) \to [0,\infty) , \epsilon > 0 \big\}$
and the collection $\big\{ c_{\epsilon}:
\epsilon > 0 \big\}$ of constants be such that, for all $\epsilon > 0$, 
$$
s_{\epsilon}(u) \sim  c_{\epsilon} u^{- \gamma}
$$
and
$$
 s_{\epsilon}(u)  \leq  s(u) \leq \big( 1 + \epsilon  \big) s_{\epsilon}(u)
$$
for all $u \geq u_0(\epsilon)$ sufficiently high. 
Then $c = \lim_{\epsilon \downarrow 0} c_{\epsilon}$ exists and
$$
 s(u) \sim c u^{-\gamma}. 
$$
\end{lemma}
\noindent{\bf Proof of Theorem \ref{theoremone}.}
Invoking Lemmas \ref{lrgthr} and \ref{lrgtwo} 
to show that the hypotheses of Lemma \ref{lrgfour}
are satisfied, we see that
$$
\lim_{k \to \infty}  \constd(k) \E_{h,\nu} \Big( \omega(T)^{\chi} \Big\vert D(T) = k \Big) \P_h  \Big( D(T) \geq k \Big)
$$
exists. Recalling further the value of the limit $\lim_{k \to \infty} \constd(k)$ from Lemma \ref{lemchipr}, the application of Lemma \ref{lrgfour} gives that

\begin{eqnarray}
&& u^{\chi}\P_{h,\nu} \Big( \omega(T) > u \Big) \nonumber\\
 &\sim& 
 \frac{1 - m_h}{\chi m_h \int_1^{\infty} y^{\chi} \log(y) d\nu(y)}
\lim_{k \to \infty} \E_{h,\nu} \Big( \omega(T)^{\chi} \Big\vert D(T) = k  \Big) \P_h  \Big( D(T) \geq k \Big)
  \nonumber \\
&=& 
 \frac{1}{\chi m_h} \Big( \int_1^{\infty} y^{\chi} \log(y) d\nu(y) \Big)^{-1}
\lim_{k \to \infty} \E_{h,\nu} \Big( \omega(T)^{\chi} 1_{D(T) = k}  \Big)
  \nonumber 
\end{eqnarray}

We thus obtain the formula (\ref{formcpo}) for the constant $\cpo$ in the statement of the theorem. $\Box$
\end{section}
\begin{section}{The proof of Theorems \ref{thmtwoaah} and \ref{theoremthree}}\label{secfive}
%We begin by proving Theorem \ref{theoremfour}. Recall that we fix a selection rule to choose a vertex $\vlast$ in the last generation $\mathcal{L}(T)$ as the starting point for the random walk. The next lemma proves that the hitting time of the root $H_{\phi}$ can be approximately coupled to an exponential random variable.
We begin by showing that the return time of a walk from the base of a typical tree to its root is well approximated by an exponential random variable. 
\begin{lemma}\label{expesc}
For any constant $\barecon > 0$, there exist $\trwtbd > 0$ and $c > 0$ such that the following holds. 
Let $T$ be a $\barecon$-bare tree satisfying $\omega(T) > \trwtbd$. 
There exists a random variable $R$ having an exponential distribution 
under  the law $\P_{T,\beta}^\vlast$, 
such that the difference
$\errzero = H_\phi - R$ satisfies, for $u >
\big(  \E_{T,\beta}^\vlast H_\phi  \big)^{1/2}$,
\begin{equation}\label{ezerbd}
\P_{T,\beta}^\vlast \Big( \big\vert \errzero \big\vert  > u \Big)
\leq \exp \Big\{ - c u^{1/40} \big( \E_{T,\beta}^\vlast H_\phi \big)^{- 1/80} \Big\}.
\end{equation}
\end{lemma}
%\noindent{\bf Proof of Theorem \ref{theoremfour}.}
%Consider the event $A$ that the tree $T$ is $\barecon$-bare. Clearly, by Lemma \ref{expesc}, on the event $A$,
%\begin{equation}\label{meanofH}
%\E_{T,\beta}^\vlast (H_\phi) = \E_{T,\beta}^\vlast (R) + O(\E_{T,\beta}^\vlast (H_\phi)^{1/2})
%\end{equation}
%Now
%
%$$ \P_{h,\nu,u} \times \P_{T,\beta}^{\vlast}  \Big(  \frac{H_\phi}{\E_{T,\beta}\big(H_\phi\big)} > t \Big) \leq 
%  \P_{h,\nu,u} \times \P_{T,\beta}^{\vlast}  \Big( \frac{H_\phi}{\E_{T,\beta}\big(H_\phi\big)} > t, A \Big) + \P_{h,\nu,u}(A^c)$$
%
%Using (\ref{meanofH}) and Lemma \ref{lemcreg}, we obtain the upper bound
%$$ 
%\limsup_{u \to \infty} \P_{h,\nu,u} \times \P_{T,\beta}^{\vlast}  \Big(  \frac{H_\phi}{\E_{T,\beta}\big(H_\phi\big)} > t \Big) \leq  \exp\{-t\}.
%$$
%Obviously the proof of the lower bound is similar. $\Box$ 

Preparing for the proof of Lemma \ref{expesc}, we decompose the duration $H_\phi$ under  $\P_{T,\beta}^\vlast$ into a large and geometrically distributed number of typically short excursions from $\vlast$, whose total duration is close to exponential, and a further brief interval of passage from $\vlast$ to $\phi$.
That is, we write the last-exit decompoosition
\begin{equation}\label{vineqn}
H_\phi = \sum_{i=1}^S N_i \, + \, F,
\end{equation}
where $S$ is the number of visits to $\vlast$ at positive times before $H_\phi$,
(so that $S$ is a geometric random variable),
and
$\sum_{i=1}^j N_i$
is the time of the $j$-th such visit, for 
$j \in \{ 1,\ldots,S \}$; and $F$
is the duration between the last visit of $X$ to $\vlast$ before time $H_\phi$
and time $H_\phi$ itself.

Before proceeding, we require
\begin{definition}\label{defnpb}
Let $\mathcal{P} = V\big(P_{\phi,\vlast}\big)$ denote the vertex-set of the 
path from the root to the base, as introduced in the decomposition in Definition \ref{defoutgr}.
 Under $\P_{T,\beta}$, let $\big\{ T_i^\mathcal{P}: i \in \N \big\}$ denote the successive distinct returns of $X$ to $\mathcal{P}$:
that is, $T_1^\mathcal{P} = \inf \big\{ i \in \N: X_i \in \mathcal{P} \big\}$, with 
$T_{n+1}^\mathcal{P} = \inf \big\{ i > T_n^\mathcal{P}: X_i \in \mathcal{P}, X_i \not=  X_{T_n^\mathcal{P}} \big\}$ for $n \geq 1$.
We abbreviate $T_i = T_i^\mathcal{P}$, and set $X_\mathcal{P}:\N \to V(T)$ by $X_\mathcal{P}(i) = X \big( T_i \big)$.
%
%We write $X_0:\{ 0,\ldots, N_1 \} \to V(T)$
%for $X$ under $\P_{T,\beta}^\vlast$ conditioned on $S \geq 1$: that is, $X_0$ records the first excursion
%of the walk from $\vlast$ in the event that this excursion is completed before $X$ first visits $\phi$.
%
%We write $\xi: \mathbb{N} \to P_{\phi,\vlast}$, $\xi = \psi \circ X$,
%for the walk $X$ projected onto $P_{\phi,\base}$.
%Similarly, we define $\xi_0: \{ 0,\ldots,N_1 \} \to P_{\phi,\vlast}$,
%$\xi_0 = \psi \circ X_0$.
\end{definition}
A few observations are useful. Recall that we write $\vert v \vert = d \big( \phi, v \big)$ for $v \in V(T)$.
\begin{lemma}\label{lemrep}
Let $T$ be a $\barecon$-bare tree.
\begin{enumerate}
\item The process $\big\{ \vert X_\mathcal{P}(i) \vert:  i \geq 0 \big\}$ 
under $\P_{T,\beta}^\vlast$ is a Markov chain. Under $\P_{T,\beta}^\vlast \big( \cdot \big\vert H_\vlast < H_\phi \big)$, the process $\big\{ \vert X_\mathcal{P}(i) \vert:  0 \leq i \leq H^{X_\mathcal{P}}_\vlast  \big\}$, until the
first return to $\vlast$ by $X_\mathcal{P}$, is a stopped Markov chain. 
\item 
There exists $c \in (0,1/2)$ such that, for all $i \in \big\{ 1,\ldots, D(T)-1\big\}$,
$$
\P_{T,\beta} \Big( \vert X_\mathcal{P}(n+1)\vert = i + 1 \Big\vert \vert X_\mathcal{P}(n)\vert = i  \Big) \geq 1/2 + c,
$$
where $n \in \N$ is arbitrary. The same statement holds under the law 
 $\P_{T,\beta}^\vlast \big( \cdot \big\vert n < H_\vlast^\mathcal{P} < H_\phi^\mathcal{P} \big)$.
%provided that $n < H^{X_P}_\vlast$.
\item there exists $c > 0$ and $\conec > 0$ such that, for any $\barecon > 0$, and for any $\barecon$-bare tree $T$,
$$
\P_{T,\beta} \Big( T_{i+1} \geq \big( \log \omega(T) \big)^{\conec \barecon} \Big\vert \vert X_\mathcal{P}(i)\vert = j \Big)
 \leq \exp \big\{ -ck \big\},
$$
for all $i \in \N$, $j \in \big\{ 1,\ldots, D(T)-1 \big\}$ and $k \in \N$. This statement also holds under 
 $\P_{T,\beta}^\vlast \big( \cdot \big\vert  n < H_\vlast^\mathcal{P} < H_\phi^\mathcal{P} \big)$.
\end{enumerate}
\end{lemma}
\noindent{\bf Proof.} The first two assertions are trivially verified.
Regarding the third, note that, if $v \in V(T) \setminus \mathcal{P}$, then the conditional distribution of $X(n+1)$
given $X(n)=v$ is the same, under the two laws 
$\P_{T,\beta}$ and  $\P_{T,\beta}^\vlast \big( \cdot \big\vert  n < H_\vlast < H_\phi \big)$. It clearly suffices then to show that there exists $\conec > 0$ such that, for any $j \in \big\{ 1,\ldots, D(T)-1 \big\}$ and $v \in V(J_i)$, the hitting time under 
$\P_{T,\beta}$ of $X$ on $\mathcal{P} \setminus \{ \psi_i \}$ given $X(0) = v$ is at most $\barecon \log \log \omega(T)$
with probability at least $\big( \log \omega(T)\big)^{-\conec \barecon}$.

To see this, 
label the vertices of the path  $P_{v,\psi_i} = \big( v = \varepsilon_0,\ldots, \varepsilon_{d(v,\psi_i)} = \psi_i \big)$.
Let $m_j$ denote the number of offspring of $\varepsilon_j$.
Note that
\begin{equation}\label{mibd}
 \sum_{j=0}^{d(v,\psi_i)} m_j \leq \vert V(J_i) \vert,
\end{equation}
since all offspring of vertices in $P_{v,\psi_i}$ belong to $V(J_i) \cup \{ \psi_{i+1} \} \setminus \{ \psi_i \}$.

Under 
$\P_{T,\beta}$, given $X(n) = v$, the walk $X$ will move along $P_{v,\psi_i}$
to $\psi_i$ in successive steps from time $n$, and then make a jump to an adjacent vertex in $\mathcal{P}$,
with probability at least 
$$
\prod_{j=0}^{d(v,\psi_i)} \frac{1}{1 + m_i \bsup}
\geq  
\prod_{j=0}^{d(v,\psi_i)} \big( 2 \bsup \big)^{-2 m_j}
 \geq   \big( 2 \bsup \big)^{-2 \vert V(J_i) \vert},
$$
where (\ref{mibd}) was used in the second inequality.
With the choice $\conec = 2 \log \big( 2 \bsup \big)$,
the statement follows, then, from $T$ being $\barecon$-bare. $\Box$ \\

The following lemma, treating the law of an excursion from $\vlast$, follows directly from Lemma \ref{lemrep}.
Its proof is left to the reader.
Recall that $N_1$ is the first return time to $\vlast$.
\begin{lemma}\label{lemnun}
Let $T$ be a $\barecon$-bare tree.
Let $\conec > 0$ be as in Lemma \ref{lemrep}(iii). There exists $c > 0$ such that, for all $v > 0$, 
$$
\mathbb{P}_{T,\beta}^\vlast \Big(  N_1 > v \Big\vert S \geq 1 \Big) \leq \exp \bigg\{ - \frac{cv}{(\log
  \omega(T))^{\conec \barecon}} \bigg\}.
$$
\end{lemma}
The following lemma has a similar proof that is also left to the reader. 
\begin{lemma}\label{lempass}
Let $T$ be a $\barecon$-bare tree.
%Under $\P_{T,\beta}$, let $H_\phi = T^{\{ \phi  \}}_1$ denote the first strictly positive visit time to $\phi$.
There exists a constant $C > 0$ such that
$$
\mathbb{E}_{T,\beta}^\phi \big( H_\vlast \big) \leq C D(T) \big( \log \omega(T) \big)^{\conec \barecon}.  
$$ 
\end{lemma}
We now bound the error in estimating the total time of excursions
from $\vlast$ by using the mean excursion time.
\begin{lemma}\label{lemell}
Let $T$ be a $\barecon$-bare tree.
Write $\mu = \E_{T,\beta}^\vlast \big( S \big)$. 
Set
$L = \sum_{i=1}^S N_i \, - \, S \mathbb{E} N_1$. 
There exists $\ctwo > 0$ such that, if 
$\omega(T) > \ctwo$ and $u \geq \mu^{1/2}$,
$$
\P_{T,\beta}^\vbase \Big( \big\vert L \big\vert \geq u  \Big)
\leq \exp \Big\{ - u^{1/40} \mu^{-1/80} \Big\}.
$$
\end{lemma}
\noindent{\bf Proof.}
It is easy to see that $\mu \geq
\binf^{D(T) - 1}$. Thus, Lemma \ref{lemom}(ii) permits us to assume that $\mu \geq
C$, where $C > 0$ is an arbitrary constant, by increasing the constant $\ctwo > 0$ as necessary.

Lemma \ref{lemnun} implies that, for some $C > 0$, under $\P_{T,\beta}^\vlast \big( \cdot \big\vert S = 1 \big)$, 
${\rm Var} \big( N_1 \big) \leq C \big( \log \omega(T) \big)^{2 \conec \barecon}$.
Hence  Theorem 3.7.1 of \cite{Dembo_Zeitouni} yields the following moderate deviations estimate.
For $0 \leq \epsilon < 1/2 $, there exists $c > 0$ such that, for all $n \in \N$,
\begin{equation}\label{stnte}
\mathbb{P}_{T,\beta}^\vlast \bigg( \Big\vert  \sum_{i=1}^n N_i - n \mathbb{E} N_1 \Big\vert 
   \geq \big( \log \omega(T) \big)^{\conec \barecon} n^{1/2 + \epsilon} \bigg\vert S = n \bigg) 
\leq \exp \big\{ - c n^{2\epsilon} \big\}.
\end{equation}
Note that, for any $k \in \N$,
\begin{eqnarray}
 & & \mathbb{P}_{T,\beta}^\vlast \Big( \vert L \vert \geq k \Big) \nonumber \\
  & = &  \sum_{i=1}^{\mu^{1 + \epsilon}} 
\mathbb{P}_{T,\beta}^\vlast \big( \vert L \vert \geq k \big\vert S = i \big) \mathbb{P} \big( S = i \big) \, + \, 
\mathbb{P}_{T,\beta}^\vlast \Big( S  > \mu^{1 + \epsilon} \Big).
  \label{llo}
\end{eqnarray}
The random variable $S$ being geometric, we have that there exists $c > 0$ such that, for any $\epsilon > 0$,
\begin{equation}\label{eqmueps}
\mathbb{P}_{T,\beta}^\vlast \Big( S  > \mu^{1 + \epsilon} \Big) \leq \exp \big\{ - c \mu^\epsilon \big\}. 
\end{equation}

For the remainder of the proof, we write $\P$ for a measure under which $\big\{ N_i: i \in \N \big\}$
is an independent and identically distributed sequence of random variables, each of which has the law of $N_1$ under
$\P_{T,\beta}^\vlast \big( \cdot \big\vert S \geq 1 \big)$.

Note that, for  
$\mu^{\frac{1}{20}} \leq i < \mu^{1 + \epsilon}$,
\begin{eqnarray}
 & & \mathbb{P}_{T,\beta}^\vlast \Big(  \vert L \vert \geq \big( \log \omega(T) \big)^{\conec \barecon}  \mu^{(1 + \epsilon)(1/2 + \epsilon)} \Big\vert S = i \Big) \nonumber \\
 & = & \mathbb{P} \Big( \Big\vert \sum_{j=1}^{i} N_j \, - \,  i \mathbb{E}
 N_1 \Big\vert \geq \big( \log \omega(T) \big)^{\conec \barecon} 
 \mu^{(1 + \epsilon)(1/2 + \epsilon)} \Big) \nonumber \\
 & \leq & \mathbb{P} \Big( \Big\vert \sum_{j=1}^{i} N_j \, - \,  i
 \mathbb{E} N_1 \Big\vert \geq \big( \log \omega(T) \big)^{\conec \barecon} 
 i^{1/2 + \epsilon} \Big) \nonumber \\
 & \leq & \exp \Big\{ - c i^{2\epsilon} \Big\}
  \leq  \exp \Big\{ - c \mu^{\frac{\epsilon}{10}} \Big\}, \label{etoi}
\end{eqnarray}
the second inequality by (\ref{stnte}) and the third by 
 $i \geq \mu^{1/20}$. Note also that,
 for $i < \mu^{\frac{1}{20}}$,
\begin{eqnarray}
 & & \mathbb{P}_{T,\beta}^\vlast \Big(  \vert L \vert \geq \big( \log \omega(T) \big)^{\conec \barecon}  \mu^{(1 + \epsilon)(1/2 + \epsilon)} \Big\vert S = i \Big) \nonumber \\
 & = & \mathbb{P} \Big( \Big\vert \sum_{j=1}^{i} N_j \, - \,  i \mathbb{E}
 N_1 \Big\vert \geq \big( \log \omega(T) \big)^{\conec \barecon} 
 \mu^{(1 + \epsilon)(1/2 + \epsilon)} \Big) \nonumber \\
 & \leq & \mathbb{P} \Big( \sum_{j=1}^{i} N_j 
 \geq 2^{-1} \big( \log \omega(T) \big)^{\conec \barecon} 
 \mu^{(1 + \epsilon)(1/2 + \epsilon)} \Big)  \nonumber \\
 & \leq &  \mu^{1/20}
 \exp \Big\{ -  c  
   \mu^{(1 + \epsilon)(1/2 + \epsilon) - 1/20} 
 \Big\}, \label{etui}
\end{eqnarray}
where, in the first inequality, we used
$$
i \mathbb{E} N_1 \leq  C \mu^{1/20} \big( \log
\omega(T) \big)^{\conec \barecon} \leq  
2^{-1} \big( \log \omega(T) \big)^{\conec \barecon}
 \mu^{(1 + \epsilon)(1/2 + \epsilon)},
$$ 
(which follows from Lemma \ref{lemnun} and the assumed $\mu \geq C$). 
The second inequality follows because one among the $N_i$ in question
exceeds the average requirement, along with Lemma \ref{lemnun}
and  $i \leq \mu^{1/20}$.
 
 By (\ref{llo}), (\ref{etoi}), (\ref{etui}) and (\ref{eqmueps}),
\begin{eqnarray}
 & & \mathbb{P}_{T,\beta}^\vlast \Big( \vert L \vert \geq \big( \log \omega(T) \big)^{\barecon \conec} \mu^{(1 + \epsilon)(1/2 + \epsilon)} \Big) \nonumber \\
 & \leq &  
    \mu^{1/10}
 \exp \Big\{ -  c  
   \mu^{(1 + \epsilon)(1/2 + \epsilon) - 1/20} \Big\}   
  \, + \, \mu^{1 + \epsilon} \exp \big\{ - c \mu^{\epsilon/10} \big\}
  + \exp \big\{ - c \mu^{\epsilon}  \big\} \nonumber \\ 
 & \leq &   \exp \big\{ - \mu^{\epsilon/20} \big\},
 \nonumber
\end{eqnarray}
the second inequality due to $\mu \geq C$.

It follows from
$\mu^{\epsilon/4} \geq \big( \log \omega(T)
\big)^{\barecon \conec}$, (which is implied by $\mu \geq \binf^{D(T)-1}$ and Lemma \ref{lemom}(ii)), that
$$
 \mathbb{P}_{T,\beta}^\vlast  \Big(  \vert L \vert \geq \mu^{1/2 + 2\epsilon} \Big)
 \leq \exp \Big\{ -  \mu^{\frac{\epsilon}{20}} \Big\}.
$$
Recalling that $u \geq \mu^{1/2}$, we may set
$\epsilon > 0 $ so that $\mu^{2\epsilon} = u \mu^{-1/2}$. In this way, we obtain the statement of the lemma. $\Box$ \\
The geometric random variable $S$ must now be compared with an exponential random variable.
\begin{lemma}\label{lexp}
Any geometric random variable may be coupled to an exponential random variable $E$
in such a way that $\vert G - E  \vert \leq 1$ almost surely.
\end{lemma}
\noindent{\bf Proof.}
To the geometric random variable $G$ 
such that $\P (G = i) = p^i(1-p)$ for $i \geq 0$, the exponential random
variable $E$ satisfying $\P \big( E > u \big) = \exp \big\{ - \log (p^{-1})
u\big\}$ for $u > 0$, may be coupled by setting $G = \lfloor E \rfloor$. $\Box$ \\

\noindent{\bf Proof of Lemma \ref{expesc}.}
By Lemma \ref{lexp}, we may construct under $\P_{T,\beta}^\vlast$ an  exponential random variable $E$
such that $\big\vert S - E \big\vert \leq 1$.
Set $R = E \mathbb{E}_{T,\beta}^\vlast \big( N_1 \big\vert S \geq 1 \big)$. 
Then $R$ has an exponential distribution, and
$$
\errzero : = \sum_{i=1}^S N_i \, - \, R
$$
satisfies
\begin{equation}\label{psstat}
\mathbb{P}_{T,\beta}^\vlast \Big( \vert \errzero \vert \geq  u  + C \big( \log \omega(T)  \big)^{\conec \barecon} \Big) \leq  \exp \Big\{ - u^{1/40} \mu^{-1/80} \Big\},
\end{equation}
by means of 
$\mathbb{E}_{T,\beta}^\vlast \big( N_1 \big\vert S \geq 1 \big) \leq C \big( \log \omega(T) \big)^{\barecon \conec}$ (which
follows from Lemma \ref{lemnun}),
 and Lemma \ref{lemell}.

Noting that $S \leq H_\phi$ on $\mathbb{P}_{T,\beta}^\vlast$, we obtain
$$
\mathbb{P}_{T,\beta} \Big( \vert \errzero \vert \geq  u  \Big) \leq  \exp \Big\{ - c u^{1/40} \big( \mathbb{E}_{T,\beta} H_\phi \big)^{-1/80} \Big\}
$$
for  $u  \geq C \big( \log \omega(T)  \big)^{\conec \barecon}$,
whence,
the statement of the lemma.
$\Box$

We now present a result that combines the content of Lemma \ref{expesc} with an asymptotic expression for  
the mean return time to the root under  $\P_{T,\beta,\falldeep}^\phi$.
The  quantity $\omegachild$ from Definition \ref{defomegachild} appears in the statement. We mention that the 
the notation of Definition \ref{defomegaext} may be employed to give a succinct expression for this quantity: $\omegachild = \sum_{v \in V(T_\vchild)} \omega_{\vchild}\big( v \big)$.
\begin{prop}\label{proptsp}
Let $\barecon \in (0,\infty)$ be an arbitrary constant. 
There exists $\trwtbd \in (0,\infty)$ and $\conthrbig > 0$, such that 
the following holds. Let $T$ denote a $\barecon$-bare weighted tree such that 
$\omegachild > \trwtbd$. 

The distribution of $H_\phi$ under $\P_{T,\beta}^{\vchild} \big( \cdot \big\vert \falldeep \big)$ is such that we may construct on this space  an exponential random variable $E$ with $\E_{T,\beta}^\vchild(E \vert \falldeep) = 2\omegachildbig/\pfd$ and with 
\begin{equation}\label{eoneerr}
 \E_{T,\beta}^{\vchild} \Big( \big\vert H_\phi - E \big\vert \Big\vert \falldeep \Big)  
 \leq \conthrbig \omegachild^{1/2}.
\end{equation}
We further have that
\begin{equation}\label{srineq}
\mathbb{E}_{T,\beta}^{\vchild} \Big( H_\phi  \Big\vert \falldeep^c  \Big)
 \leq C \Big( \log \omegachild \Big)^{\barecon \conec}.
\end{equation}
\end{prop}
\noindent{\bf Proof.}
We begin by showing (\ref{srineq}). Note that, under $\P_{T,\beta}^\vchild \big( \cdot \big\vert \falldeep^c \big)$,
the process $\vert X_P \vert$ (which is specified in Definition \ref{defnpb}), 
has the conditional distribution of 
 $\vert X_P \vert$  under $\P_{T,\beta}^\vchild$ given that this process hits zero before $D(T)$. The law 
of $\vert X_P \vert$  under $\P_{T,\beta}^\vchild$ having a uniform bias to the right (except for a reflection at $D(T)$), we see that this conditional distribution has a uniform bias to the left, so that the hitting time of $0$ by this process has a finite mean, uniformly in $T$. Note further that the bound in 
Lemma \ref{lemrep}(iii) remains valid under the law  $\P_{T,\beta}^\vchild \big( \cdot \big\vert \falldeep^c \big)$,
since the walk $X$ under this measure has the unconditioned jump distribution at points in $V(T) \setminus P$. In this way, we conclude that (\ref{srineq}) holds.

To prove (\ref{eoneerr}), we now argue that
\begin{equation}\label{qrineq}
 \mathbb{E}_{T,\beta}^{\vchild} \Big( H_\vlast  \Big\vert \falldeep \Big) \leq C
 \Big(\log \omegachildbig \Big)^{\conec \barecon + 1}.
\end{equation} 
Indeed, using the notation (\ref{vineqn}), under $\P_{T,\beta}^\vchild$, 
we have that $\falldeep = \big\{ S \geq 1 \big\}$. By Lemma \ref{lempass}, 
$E_{T,\beta}^\vchild \big( H_\vlast \big) \leq C D(T) \big(\log \omegachildbig \big)^{\conec \barecon}$.
We now obtain (\ref{qrineq}) from $\binf^{D(T_0)-1} \leq \omegachild$.

By (\ref{eqelres}), we have that
\begin{eqnarray}
 & & 2 \omegachildbig - 1  =  \mathbb{E}^{\vchild}_{T,\beta} \big( H_\phi \big) \nonumber \\
 & = &
\mathbb{E}_{T,\beta}^{\vchild} \Big( H_\phi  \Big\vert \falldeep \Big) \mathbb{P}^{\vchild}_{T,\beta} \Big( \falldeep \Big) 
%\nonumber \\
% & & \qquad \qquad \qquad + \, 
+ \mathbb{E}_{T,\beta}^{\vchild} \Big( H_\phi  \Big\vert \falldeep^c \Big)
\mathbb{P}_{T,\beta}^{\vchild} \Big(  \falldeep^c  \Big) \nonumber \\
& = & \pfd
\bigg( \mathbb{E}_{T,\beta}^{\vchild} \Big( H_\phi - H_\vlast \Big\vert \falldeep \Big)
 + \mathbb{E}_{T,\beta}^{\vchild} \Big( H_\vlast  \Big\vert \falldeep \Big) \bigg) \nonumber \\
   & & \qquad \qquad 
 \, + \, \big(1-\pfd\big) 
\mathbb{E}^{\vchild}_{T,\beta} \Big( H_\phi   \Big\vert \falldeep^c  \Big). \nonumber 
\end{eqnarray}
Hence,
\begin{eqnarray}
& & \bigg\vert   
\mathbb{E}_{T,\beta}^{\vchild} \Big( H_\phi - H_\vlast \Big\vert \falldeep \Big)  \, - \,
\frac{2\omegachildbig - 1}{\pfd} \bigg\vert \label{ecsq} \\
 & \leq & 
\mathbb{E}_{T,\beta}^{\vchild} \Big( H_\vlast  \Big\vert \falldeep \Big) \, + \, \big( \pfd^{-1} - 1 \big)
\mathbb{E}_{T,\beta}^{\vchild} \Big( H_\phi  \Big\vert \falldeep^c   \Big) 
  \leq   C
 \big( \log \omegachild \big)^{\barecon \conec + 1}, \nonumber
\end{eqnarray}
where, in the second inequality, we used (\ref{qrineq}) and (\ref{srineq}), as well as  that $\pfd$  is bounded away from zero uniformly, (as shown by the readily verified inequality, $\pfd \geq 1 - \binf^{-1}$).

By Lemma \ref{expesc}, under $\P_{T,\beta}^\vchild \big( \cdot \big\vert \falldeep \big)$,
we may construct an exponentially distributed random variable $\reo$ 
in such a way that
$H_\phi = H_\vlast + \reo + \ezero$, where
$\ezero$ satisfies the bound 
\begin{equation}\label{eqeo}
\E_{T,\beta}^\vchild \Big( \big\vert \ezero \big\vert \Big\vert \falldeep \Big) \leq C \E_{T,\beta}^\vlast \big( H_\phi \big)^{1/2}.
\end{equation}
%\ref{ezerbd}). 
By (\ref{ecsq}) and (\ref{eqeo}),
\begin{eqnarray}
  \Big\vert  \E_{T,\beta}^\vchild  \big( \reo \big\vert \falldeep \big)  - 2\pfd^{-1} \omegachild \Big\vert 
 & \leq & C \big( \log \omegachild \big)^{\barecon \conec + 1} + \E_{T,\beta}^\vchild \Big( \vert \ezero \vert \Big\vert \falldeep \Big)
\nonumber \\
 & \leq & C \big( \log \omegachild \big)^{\barecon \conec + 1} + C \big( \E_{T,\beta}^\vlast H_\phi
 \big)^{1/2}.
  \label{evomq}
\end{eqnarray}
From $\pfd \E_{T,\beta}^\vlast (H_\phi) \leq \E_{T,\beta}^\vchild (H_\phi)$, $\pfd \geq 1 - \binf^{-1}$ and (\ref{eqelres}), we see that 
\begin{equation}\label{eonetwo} 
\E_{T,\beta}^\vlast (H_\phi) \leq C \omegachildbig.
\end{equation}
This implies that
\begin{equation}\label{evtw}
\Big\vert  \E_{T,\beta}^\vlast \big( \reo \big\vert \falldeep \big)  - 2 \pfd^{-1} \omegachildbig \Big\vert
 \leq C \omegachildbig^{1/2},
\end{equation}
since $\omegachildbig$ is assumed to exceed a large constant $\trwtbd$.
Now, we let $$
\reotilde = \frac{2 \omegachildbig}{\pfd} \frac{\reo}{\E_{T,\beta}^\vlast \big( \reo \big\vert \falldeep \big) },
$$
so that $\reotilde$ under $\P_{T,\beta}^\vchild \big( \cdot \big\vert \falldeep \big)$ is an exponential random variable with 
$\E_{T,\beta}^\vchild \big( \reotilde \big\vert \falldeep \big) =  2 \pfd^{-1}\omegachildbig$. Thus, under  $\falldeep$, we have that 
\begin{equation}\label{hphirep}
H_\phi = \reotilde + \eone,
\end{equation}
where $\eone = H_\vlast +  \big( \reo - \reotilde \big) + \ezero$ satisfies
$$ 
\E^{\vchild}_{T,\beta} \Big( \big\vert \eone \big\vert \Big\vert \falldeep \Big)
 \leq
 C \big( \log \omegachildbig
 \big)^{\barecon \conec + 1}
+
 C \omegachildbig^{1/2} 
    + C \omegachildbig^{1/2} 
 \leq C \omegachildbig^{1/2},
$$
the first inequality by 
(\ref{qrineq}), (\ref{evtw}), (\ref{eqeo}) and (\ref{eonetwo}), and the second by $\omegachildbig > \trwtbd$. 
We have obtained (\ref{eoneerr}). $\Box$ \\
\noindent{\bf Proof of Theorems \ref{thmtwoaah} and \ref{theoremthree}.}
Note that Theorem \ref{thmtwoaah} follows directly from the two statements in Theorem \ref{theoremthree}.

Note that the distribution of $H_\phi$ under $\P_{T,\beta}^\vbase$
coincides with that of $H_\phi - H_{\vbase}$ under $\P_{T,\beta,\falldeep}^\vchild$.
For this reason, the first statement of Theorem \ref{theoremthree} follows from (\ref{ecsq}) and Lemma \ref{lemcreg}.

We begin deriving the second statement of Theorem \ref{theoremthree}. 
Let $T$ be any weighted tree. Clearly, under the conditional measure 
$\P_{T,\beta,\falldeep}^\phi$, we have that $X_1 = \vchild$ almost surely, so that the distribution of $H_\phi$ under $\P_{T,\beta,\falldeep}^\phi$
coincides with that of $1 + H_\phi$ under $\P_{T,\beta,\falldeep}^\vchild$. For the second statement of  
Theorem \ref{theoremthree}, it thus suffices to show that
\begin{equation}\label{eqvchild}
 \Big( \P_{h,\nu,u} \times \P_{T,\beta,\falldeep}^\vchild \Big)
  \Big( \frac{H_\phi}{2\omegachild \pfd^{-1}} > t \Big) \to \exp \big\{ -t \big\}
\end{equation}
as $u \to \infty$. 
To demonstrate this, let $A$ denote the event that $T$ is $\barecon$-bare and that $\omegachild > u^{\frac{\log \binf}{4 \log \bsup}}$. Note that, by $\omegachild \geq \binf^{D(T) - 1}$ and Lemma \ref{lemom}, we have that a $\barecon$-bare tree $T$ for which $\omega(T) \geq u$ satisfies  $\omegachild >  
u^{\frac{\log \binf}{4 \log \bsup}}$. By this and Lemma \ref{lemcreg}, we find that
\begin{equation}\label{eqproba}
 \P_{h,\nu,u}\big( A \big) \to 1
\end{equation}
as $u \to \infty$. Under the law $\P_{h,\nu,u} \times 
\P_{T,\beta,\falldeep}^\vchild \big( \cdot \big\vert A \big)$, we construct a random variable $E$ that, conditioning on the tree, enjoys the properties given in Proposition \ref{proptsp}. Fixing $\epsilon > 0$, we have that 
\begin{equation}\label{eqthreea}
 \Big( \P_{h,\nu,u} \times \P_{T,\beta,\falldeep}^\vchild \Big)
  \Big( H_\phi  > 2 \omegachild \pfd^{-1} t \Big) \leq A_1 + A_2 + A_3,
\end{equation}
where
$$
A_1  =  \Big( \P_{h,\nu,u} \times \P_{T,\beta,\falldeep}^\vchild \Big)
  \Big( E  > 2 \omegachild \pfd^{-1} t(1-\epsilon), A  \Big),
$$
$$
A_2 =  \Big( \P_{h,\nu,u} \times \P_{T,\beta,\falldeep}^\vchild \Big)
  \Big( \big\vert H_\phi - E \big\vert  > 2 \omegachild \pfd^{-1} t\epsilon, A  \Big)
$$
and
$A_3 = \P_{h,\nu,u}  \big( A^c \big)$.
Note that (\ref{eqproba}) says that
$A_3 \to 0$ as $u \to \infty$. Recalling that, conditionally on $T$, $E$ is an exponential random variable of mean $2\omegachild/\pfd$, (\ref{eqproba}) yields that
$$
 \lim_{u \to \infty} A_1 = \exp \big\{ - t(1-\epsilon) \big\}.
$$ 
To bound $A_2$, fix a  $\barecon$-bare tree for which $\omegachild >  u^{\frac{\log \binf}{4 \log \bsup}}$.
By Proposition \ref{proptsp} and Markov's inequality, we have that
$$
\P_{T,\beta,\falldeep}^\vchild \Big(  \big\vert H_\phi - E \big\vert \geq  2 \omegachild \pfd^{-1} \epsilon t  \Big)
 \leq \conthrbig \big( \epsilon t \big)^{-1} \omegachild^{-1/2}
  \leq  \conthrbig \big( \epsilon t \big)^{-1} u^{- \frac{\log \binf}{8 \log \bsup}}.
$$
Thus,
$A_2 \leq \conthrbig \big( \epsilon t \big)^{-1} u^{- \frac{\log \binf}{8 \log \bsup}}$.
Substituting the obtained bounds into (\ref{eqthreea}), and taking limits $u \to \infty$ followed by $\epsilon \to 0$, yields
$$
\limsup_{u \to \infty} \Big( \P_{h,\nu,u} \times \P_{T,\beta,\falldeep}^\vchild \Big)
  \Big( H_\phi  > 2 \omegachild \pfd^{-1} t \Big) \leq \exp \big\{ -  t \big\}. 
$$
The complementary lower bound has a verbatim proof. In this way, we obtain (\ref{eqvchild}). This completes the proof of Theorem \ref{theoremthree}. $\Box$
\end{section}
\appendixpage
\appendix
\begin{section}{Coding of trees}\label{appa}
%In Appendix \ref{appb}, we will prove Theorem \ref{lemdivexpold}. The proof will require some variations on the framework used to define finite rooted trees and the Galton-Watson distribution on them. 

We provide here a precise formulation of the notion of a finite rooted tree and the Galton-Watson distribution on them. This treatment has been taken essentially from \cite{LeGall}. 

We begin by defining a space of labels
$$
\mathcal{U} = \bigcup_{n=0}^\infty \N^n,
$$
where $\N = \big\{ 1,2,\ldots \big\}$ and, by convention, $\N^0 = \{ \emptyset \}$. Each element $u = (u_1,\ldots,u_n) \in \mathcal{U}$ is thus a finite sequence of natural numbers. We denote by $\vert u \vert = n$ the generation of $u$. We define the concatenation $uv$ of $u = (u_1,\ldots,u_n)$ and $v = (v_1,\ldots,v_m)$ by means of $uv = (u_1,\ldots,u_n,v_1,\ldots,v_m)$.

The mapping $\pi:\mathcal{U} \setminus \{ \emptyset \} \to \mathcal{U}$, given by $\pi \big( u_1\ldots u_n \big) = u_1 \ldots u_{n-1}$ associates to each individual its parent.

A finite rooted ordered tree $T$ has vertex set given by a finite subset $\mathcal{U}$ such that
\begin{enumerate}
\item $\emptyset \in T$.
\item $u \in T \setminus \{ \emptyset \}$ implies that $\pi(u) \in T$.
\item For every $u \in T$, there exists an integer $k_u(T) \geq 0$ such that, for every $j \in \N$, $uj \in T$ if and only if $1 \leq j \leq k_u(T)$.
\end{enumerate}
The number $k_u(t)$ is interpreted as the number of offspring of $u$ in $T$.

The edges of such a tree are the (unoriented) edges connecting each of its elements $u$ (except $\emptyset$) to its parent $\pi(u)$. Such a tree $T$ is weighted when we associate to it a function $\beta:E(T) \to (0,\infty)$.

The formal definition of the Galton-Watson law is now given. 
\begin{definition}
Let $h = \big\{ h_i: i \in \N \big\}$ satsfy $\sum_{i=1}^\infty i h_i \leq 1$. 
Define a family of independent 
and identically distributed 
random variables $\big\{ K_u: u \in \mathcal{U} \big\}$, each having the law $h$. Write $\theta$ for the subset of $\mathcal{U}$ given by
$$
\theta = \Big\{ u = \big( u_1\ldots u_n \big) \in \mathcal{U}: 
 u_j \leq K_{u_1\ldots u_{j-1}} \, \textrm{for each $1 \leq j \leq n$} \Big\}.
$$
Then $\theta$ is the vertex set of the Galton-Watson tree sampled according to $\P_h$. 
%The edges of the tree are the (unoriented) edges $\big( u,\pi(u) \big)$ for $u \in \theta \setminus \{ \emptyset \}$. 

With $\bsup \geq \binf > 1$, and for a law $\nu$ supported in $[\binf,\bsup]$, 
the law $\P_{h,\nu}$ is defined by introducing an independent collection $\big\{ \beta_u: u \in \mathcal{U} \setminus \{ \emptyset \} \big\}$ of $\nu$-distributed random variables, and setting  $\beta_e = \beta_u$, for each edge $e = (u,\pi(u))$, with $u \in \theta \setminus \{ \emptyset \}$.
\end{definition}
Note that $\P_{h,\nu}$ is a law on finite rooted weighted ordered trees.
In the next appendix, we will wish to work with the unordered variant of this object. We now formally define this.
\begin{definition}\label{defformphnu}
Let $T$ and $T'$ be two finite rooted weighted ordered trees. We say that $T$ and $T'$ are isomorphic if there exists a bijection $\psi:V(T) \to V(T')$ with the properties that
\begin{enumerate}
\item $\psi(\emptyset) = \emptyset$.
\item $\psi \circ \pi = \pi \circ \psi$ on $V(T)$.
\item $\beta(\psi(e)) = \beta(e)$ for each $e \in E(T)$.
\end{enumerate}
In 3, we used 2 to extend the definition of $\psi$ to $E(T)$: indeed,
each edge $e \in E(T)$ taking the form $e = \big( \pi(u) , u \big)$ for some 
$u \in V(T) \setminus \{ \emptyset \}$, we set
$\psi(e)$ to be the edge $\big( \psi \circ \pi (u),\psi(u) \big)$.

A rooted weighted {\it unordered} tree is an isomorphism class of such trees. 
An unordered weighted Galton-Watson tree is the law of 
the isomorphism class of a sample of $\P_{h,\nu}$.   
\end{definition}

After Definition \ref{defphnu}, we introduced the joint tree and walk measure  $\P_{h,\nu} \times \P_{T,\beta}^v$. The definition was left imprecise, in the sense that a {\it selection rule} was required to specify the vertex $v \in V(T)$. We are now able to clarify this point formally.
\begin{definition}
Let $\Theta$ denote the set of finite rooted weighted ordered trees. A selection rule is a map $s: \Theta \to \mathcal{U}$ with the property that $s(T) \in V(T)$ for all $T \in \Theta$. In the definition of  $\P_{h,\nu} \times \P_{T,\beta}^v$, $v$ is a selection rule.  
\end{definition}
\end{section}
\begin{section}{The proof of Theorem \ref{lemdivexpold}}\label{appb}
We aim to show that, for any $u > 0$, under the law $\P_{h,\nu,u}$, (which we recall denotes $\P_{h,\nu}\big( \cdot \big\vert \omega(T) > u \big)$), a long outgrowth is unlikely. 
\begin{subsection}{Surgery using the FSO-decomposition}
The plan is to argue that, under the measure $\P_{h,\nu}$, such an outgrowth is not the most efficient means of securing the condition $\omega(T) > u$ of high weight. We will demonstrate that, 
from a high-weight tree with a long outgrowth, 
the outgrowth may be removed, 
and the tree lengthened a little, 
in such a way that the surgically altered tree 
has at least the weight of the original tree, 
with the altered tree being demonstrably more probable 
under $\P_{h,\nu}$ than the original one.

The proof that we are explaining being a little involved, 
we prefer to give it firstly under two assumptions, 
that serve to remove some distracting details in the argument.
\begin{hyp}\label{hyplem}
The edge-weight law $\nu$ has no atoms: that is, for all $x \in (1,\infty)$, $\nu \big( \{ x \} \big) = 0$.
The offspring distribution $h = \big\{ h_i: i \in \N \big\}$ satisfies $h_1 > 0$. 
\end{hyp}
The first hypothesis is being used because it ensures that the weights attached to the vertices in any finite tree are distinct. The second allows the presence in a tree of finite paths each of whose vertices has a single offspring. Regarding hypotheses, we also mention that we will make the harmless assumption that $h_i > 0$ for some $i \geq 2$: for, were this to fail, all of our trees would be finite paths, so that Theorem \ref{lemdivexpold}.

To define and analyse the surgical procedure, we introduce a new decomposition of a tree.
\begin{definition}\label{defdecomplem}
Let $T$ denote a weighted tree such that the maximal value of $\omega(v)$ among $v \in V(T)$ is assumed by a unique vertex that we will denote by $\vmax$. Recall that $\phi = \phi_T$ denotes the root of $T$.
 
Let the first branch point $\fbp$ denote the first vertex on the path $P_{\phi,\vmax}$ from $\phi$
to $\vmax$ having at least two offspring. (We set $\fbp = \vmax$ if there is no such vertex in $P_{\phi},\vmax$.) We define the {\it foundation}  $F$ to be $P_{\phi,\fbp}$. We define the spine $S$ to be $P_{\fbp,\vmax}$.

Writing $s = \big\vert V(S) \big\vert - 1$, we further label the successive vertices in $S$ as
$\fbp = \chi_0,\chi_1,\ldots,\chi_s = \vmax$. For $i \in \big\{ 0,\ldots, s \big\}$, we let $O_i$, the $i$-th offshoot of $T$, denote the connected component containing $\chi_i$ in the graph with vertex set $V(T)$ and edge-set $E(T) \setminus E(S)$. We also set $O_i = \emptyset$
 if $i > s$. (Note that $O_s$ is equal to the singleton graph with vertex $\vmax$.)

In this way, the edges of $T$ are partitioned into the foundation $F$, the spine $S$ and the various offshoots $O_j$, $0 \leq j \leq s - 1$. We refer to the decomposition as the $\fso$-decomposition. 
\end{definition}
\begin{figure}\label{figdecomplem}
\begin{center}
\includegraphics[width=0.2\textwidth]{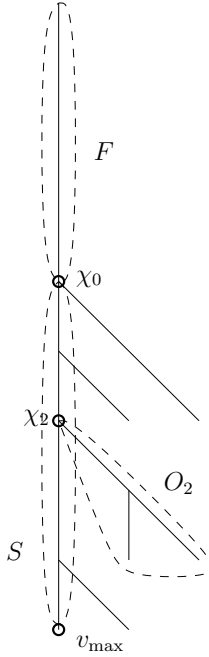} \\
\end{center}
\caption{The $\fso$-decomposition 
%given in Definition \ref{defdecomplem} 
is illustrated.}
\end{figure}
To obtain Theorem \ref{lemdivexpold}, our plan is to show that any given offshoot $O_i$ is likely to be small under $\P_{h,\nu,u}$. 
A short additional argument will then yield the same conclusion for any given outgrowth $J_i$, as desired. (The offshoots and the outgrowths are each the subtrees hanging off the path connecting the root to a vertex near the base of the tree, $\vmax$ or $\vbase$. Except for a few terms at the extreme ends, the list of outgrowths coincides with that of the offshoots up to a small additive random shift in the index. In other words, the offshoots are a convenient technical tool, and it is only a small step to learn that the outgrowths are small once we know that the offshoots are.) 

To show that the offshoots are small, we will fix $i \in \N$, and condition a sample $T$ of $\P_{h,\nu}$ on all elements in its $\fso$-decomposition except for $F$ and $O_i$, and then show that, for this law, there is a more probable means to achieve $\omega(T) > u$ than by insisting that $O_i$ be large. This means is to insist instead that $F$ be long (but with a length that is much shorter than the size of $O_i$ under the comparison, so that this outcome is the more probable).

\begin{figure}\label{figsurgery}
\begin{center}
\includegraphics[width=0.5\textwidth]{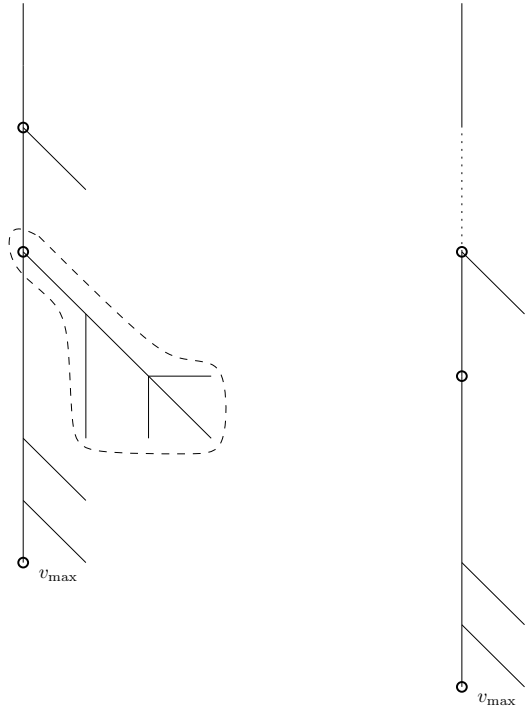} \\
\end{center}
\caption{Before and after: a large offshoot is removed. The only trace that the procedure leaves is the unobtrusive short dotted segment shown on the right.}
\end{figure}

Figure 2
%\ref{figsurgery} 
illustrates the surgical procedure that is the main tool in the proof of Theorem \ref{lemdivexpold}. The left-hand figure depicts a tree sampled under $\P_{h,\nu,u}$ that happens to contain a large offshoot. We now describe the surgery in a little more detail: the formal argument corresponding to this sketch will appear in the proof of Lemma \ref{lemofsc}. We will associate to the tree on the left-hand-side, another one, shown on the right, that still satisfies the condition that the tree weight be at least $u$, but which is more probable under $\P_{h,\nu}$ than is the original tree.
This latter tree will be obtained from the first by some surgery. 
The offending offshoot will be trimmed from the tree. 
This trimming necessarily entails some loss in the weight of the tree. 
This weight must be regained, if the outcome is to have weight at least $u$. 
The means of gaining weight is to add a few edges into the foundation of the tree. 
The number of edges that we need to add is only logarithmic in the size of the removed offshoot, so that the cost of this alteration to the tree is negligible compared to the cost of the original offshoot. It is in order to ensure that only a short lengthening of the foundation is sufficient to regain the lost weight that we have chosen to use the FSO-decomposition, defined in terms of the maximum weight vertex, $\vmax$. Indeed, the removed offshoot necessarily does not contain $\vmax$, so that, if this offshoot has a certain size $k$, the weight of the tree formed by its removal  is at least $1/(k+1)$ that of the original weight. The fact that this trimmed tree necessarily has a significant weight is then vital in arguing that a slight extension to the foundation of this tree yields one with at least the weight of the original tree.

This , then, is the plan in outline. To implement it, we need a little more notation.
Firstly, recall that,
in Appendix \ref{appa}, we defined the weighted Galton-Watson law $\P_{h,\nu}$ as a measure on 
ordered weighted trees. As we will shortly explain, the ordering information specified by this definition is inconvenient for our present method of the proof. As such, we wish to consider unordered trees instead. We henceforth abuse notation and write $\P_{h,\nu}$ for the law on unordered weighted trees given in Definition \ref{defformphnu}. 
\begin{definition}\label{defphnui}
 Let $s \in \N$ and $i \in \big\{0,\ldots,s-1\big\}$. Let $S'$  denote a weighted path of length $\vert E(S') \vert =s$, and let $O'_j$, $j \in \big\{ 0,\ldots,s-1\big\} \setminus \{ i \}$ denote a collection of (rooted) weighted trees, for which there exists a (rooted) weighted tree $T$ whose $\fso$-decomposition satisfies $S = S'$, and $O_j = O'_j$, for such $j$. We will refer to the collection of data
$\big(s,i,S', \{ O_j' :j \not= i \} \big)$ as an $i$-absent foliage. For brevity, we write
$\phni$ for the conditional distribution of $\P_{h,\nu}$ given that the $\fso$-decomposition of $T$
satisfies  $S = S'$ and $O_j = O'_j$ for such $j$. Note that implicit in the notation $\phni$ is a given $i$-absent foliage.  
\end{definition}
The deduction of our use of surgery, that rules out long offshoots in the $\fso$-decomposition, is now stated.
The remainder of this subsection is devoted to establishing the result.
\begin{prop}\label{lemsmallos}
Assume Hypothesis \ref{hyph}, and that $\nu$ has support in $[\binf,\bsup] \subseteq (1,\infty)$.
There exists $c > 0$ such that, for any $u > 0$,
for all $i,k \in \N$, and for all $i$-absent foliages,
$$
\phni \Big( \big\vert V(O_i) \big\vert = k \Big\vert  \omega(T)> u \Big) \leq \exp \big\{ - ck \big\}.
$$ 
\end{prop}
\noindent{\bf Remark.} From the form of this lemma, we may explain why we work with unordered trees. The statement would be untrue in the other case. For an 
ordered tree, each edge in the spine of the tree's $\fso$-decomposition would carry with it a positive integer indicating its order among the edges emanating from its parent vertex. If this number were high, then all the edges with a lower index would necessarily live in a single offshoot. 
 This would mean that certain offshoots would necessarily be large, so that the statement of Proposition \ref{lemsmallos} would be false for ordered trees.

To prove Proposition \ref{lemsmallos}, we fix $k \in \N$, and aim to show that, for a given $i$-absent foliage, 
it is probabilistically cheaper to obtain the condition that $\omega(T) > u$ by sampling a slightly long foundation  (of length about $\log k$) than by sampling a very large offshoot $O_i$ (of size $k$).
In seeking to realize this aim, note the following difficulty. If $u > 0$ is large, and the spine and offshoots fixed in the $i$-absent foliage are small, then it will be necessary for $\omega(T) > u$ that $F$ be long (having a length of order $\log u$) if the offshoot $O_i$ has a size of at most $k$. In considering a lengthening of $F$ as an alternative means to $\vert V(O_i) \vert = k$ for attaining $\omega(T) > u$, we are referring to the additional length of $F$ in excess of this necessary length. 

We now prepare to give a precise definition of a splitting of $F$ into an initial ``necessary'' subpath and a final ``optional'' subpath. 
\begin{definition}\label{defnot}
Note that, under the conditional law $\phni$, the tree $T$ is a random variable specified by the data $F$
and $O_i$. We will sometimes denote this dependence explicitly: $T = T\big( F,O_i \big)$. Note further that, under this law, the descendent tree $T_{\fbp}$ is a random variable specified by $O_i$. We will sometimes write $T_\fbp = T_\fbp \big( O_i \big)$.
\end{definition}
\begin{definition}
 Given a weighted tree $S$ and a subset $A \subseteq V(S)$, we set $\omega^S(A) = \sum_{v \in V(S)}\omega_{\phi(S)}(v)$, where the summand is written in the notation of Definition \ref{defomegaext}. We will make choices of $S$ and $A$ in terms of the notation in Definition \ref{defnot} that specifies 
tree-valued random variables under the law $\phni$. For example, given the data implicit in 
$\phni$, the quantity  $\omega^{T_{\fbp}(\emptyset)}\Big( T_{\fbp}\big(\emptyset\big) \Big)$
denotes $\sum_{u \in V(T_{\fbp})}\omega_{\fbp}(u)$ when the choice $O_i = \emptyset$ is made. That is,
$$
\omega^{T_{\fbp}(\emptyset)}\Big( T_{\fbp}\big(\emptyset\big) \Big)
 = \sum_{j=0}^s \omega_{\fbp}(\chi_j) \, + \,  \sum_{j=0,j \not=i}^s \sum_{v \in V(O_j),v \not= \chi_j}
 \omega_{\fbp}(v).
$$
\end{definition}
%In Proposition \ref{lemsmallos}, we condition $\P_{h,\nu}$ on $\omega(T > u)$, and on all parts of the decomposition, except for the $i$-th offshoot $O_i$ and the foundation $F$. To explain the idea of the proof of the lemma, suppose for a moment that we condition $\P_{h,\nu}$ on the same set of information, excepting the requirement $\omega(T) > u$. If we then consider how the condition $\omega(T) > u$ might be met under the conditonal law, we have two alternatives: either $O_i$ is large, or $F$ is long. We will argue that the latter alternative is probabilistically cheaper. Considering the event $\big\vert O_i \big\vert = k$ to mean that $O_i$ is large, we will obtain Proposition \ref{lemsmallos}.

%\begin{lemma}\label{lemnecf}
%Let $i \in \N$, $s \geq i$, $S = \big[ \chi_0,\ldots,\chi_s \big]$ and 
%$\O_j$, for $j \in \big\{ 0,\ldots,s \big\} \setminus \{ i \}$, be given.
%Let $T$ be a rooted weighted tree whose decomposition has this data, and for which $\omega(T) > u$.
%\end{lemma}
At present, we seek to prove Proposition \ref{lemsmallos} under the assumption of Hypothesis \ref{hyplem}. In the ensuing lemmas, we assume this hypothesis, as well as the other hypotheses of Proposition \ref{lemsmallos}. 
When we are done, we will revisit these arguments, to remove the need for Hypothesis \ref{hyplem}.
\begin{lemma}\label{lemnecf}
Let $T$ be a weighted tree whose $\fso$-decomposition has 
a given $i$-absent foliage $S'$, $O_j'$, $j \in \big\{ 0,\ldots,s-1\big\} \setminus \{ i \}$ for $s > i \geq 0$. 
Suppose further that  $\omega(T) > u$ and that $\vert V(O_i) \vert \leq k$. Then $T$ has the property that the inequality
\begin{equation}\label{impineq}
   \frac{\binf}{\binf - 1} \bigg( \omega^{T_{\fbp}(\emptyset)}\Big( T_{\fbp}\big(\emptyset\big) \Big)  \, + \, \omega_{\fbp}(\chi_i) \sum_{\ell=1}^k \bsup^\ell \bigg) \omega \big( v \big) > u.
\end{equation}
is satisfied by the choice $v = \fbp$.
Note that $T_{\fbp}\big(\emptyset\big)$ denotes the subtree of $T$ induced by the union of the vertices of $S$ and of the offshoots excluding $O_i$.
\end{lemma}
The following definition is convenient.
\begin{definition}\label{deflkellk}
For $k \in \N$, let $L_k$ denote the weighted tree consisting of $k$ consecutive edges, each having bias $\bsup$.
Let $\ell_k$ denote the weighted tree  consisting of $k$ consecutive edges, each having bias $\binf$.
\end{definition}
\noindent{\bf Proof of Lemma \ref{lemnecf}.}
Note firstly that, for such a tree $T$, 
$$
 \omega(T) = \sum_{v \in V(F), v \not= \fbp} \omega(v) \, + \, \sum_{v \in V(T_\fbp)} \omega(v).
$$
Note that, for $u \in V(F)$, $\omega(u) \leq \binf^{- d(u,\fbp)}
   \omega \big( \fbp \big)$, since all edge have bias are at least $\binf$. Hence,
$$
\sum_{v \in V(F), v \not= \fbp} \omega(v) \leq \frac{1}{\binf - 1} \omega (\fbp)
 \leq \frac{1}{\binf - 1}  \sum_{v \in V(T_{\fbp})} \omega (v),
$$
whence
$$
\omega(T) \leq \Big( 1 + \frac{1}{\binf - 1} \Big)  \sum_{v \in V(T_{\fbp})} \omega (v).  
$$
For $v \in V(T_{\fbp})$, $\omega(v) = \omega(\fbp) \omega^{T_{\fbp}}(v)$,
implying that
\begin{equation}\label{vkineq}
\omega(T) \leq \frac{\binf}{\binf - 1} \omega \big( \fbp \big) \omega^{T_{\fbp}}\big( T_{\fbp} \big).
\end{equation}
Note that the tree $L_k$ is the $(k + 1)$-vertex tree of maximal weight.
Adopting the notation of Definition \ref{defnot}, note thus that, among choices of $O_i = O$ for which 
$\vert V(O) \vert \leq k$, 
we have that
\begin{equation}\label{vvineq}
 \omega^{T_{\fbp}(O)}\Big( T_{\fbp}\big(O\big) \Big)
    \leq  \omega^{T_{\fbp}(\emptyset)}\Big( T_{\fbp}\big(\emptyset\big) \Big)
  +  \omega^{T_{\fbp}(\emptyset)}(\chi_i) \sum_{\ell = 1}^k \bsup^\ell.  
\end{equation}
From (\ref{vkineq}), (\ref{vvineq}) and the assumption that
$\big\vert V\big(O_i\big) \big\vert \leq k$, we obtain the statement of the lemma. $\Box$ \\
Lemma \ref{lemnecf} permits the following definition.
\begin{definition}
Let $T$ be a weighted tree satisfying the hypotheses of Lemma \ref{lemnecf}. Let $\vnec$
denote the vertex $v$ in $F$ closest to $\phi$ among those satisfying (\ref{impineq}).
We refer to $\vnec$ as the essential vertex. We call the path $P_{\phi,\vnec} \subseteq  F$, the essential foundation $\ef$. We call the path $P_{\vnec,\fbp} \subseteq F$, the optional foundation $\opf$. 
\end{definition}
We are aiming to prove Proposition \ref{lemsmallos}, which concerns  the conditional law $\phni$, in which we condition on a given $i$-absent foliage. To do so, it turns out to be better to condition also on the essential foundation, since what is then left undetermined is the $i$-th offshoot and the optional foundation, and it is in terms of the length of the optional foundation that we will shortly phrase a convenient sufficient condition for $\omega(T) > u$. We now extend the $\phni$ notation to incorporate this additional conditioning.
\begin{definition}\label{defphnuiplus}
Suppose given an $i$-absent foliage. Let $u > 0$. Under the law $\phni$,
we write $\ess$ for the event that $\vnec$ exists. (Note that this event depends implicitly on the $i$-absent foliage, on $k \in \N$ and on $u$.) Let $\ef'$ denote a fixed weighted tree (which is necessarily has the form of a succession of edges) which is a possible choice for $\ef$  under the law $\phni \big( \cdot \big\vert \ess \big)$.
We call the data consisting of the given $i$-absent foliage, and $\ef'$,
 an $i$-absent verdant foliage.
We adopt the notation $\phniplus = \phni \big(\cdot \big\vert \ess, \ef = \ef' \big)$. Note that implicit in the definition of $\phniplus$ is an $i$-absent verdant foliage, which itself depends implicitly on $k \in \N$ and $u > 0$.
\end{definition}
To summarise, for a sample $T$ of $\P_{h,\nu}$, given the spine $S$ and the offshoots $O_j$, and with $O_i$ unspecified except for the condition that $\big\vert V(O_i) \big\vert \leq k$, it is required, for $\omega(T) > u$, that a certain initial segment of the foundation $F$ exist. This segment is $\ef$. What is left of $F$ after $\ef$ is $\opf$. By further conditioning $T$ on $\ef$, (that is, by using the law $\phniplus$), we have arrived at a convenient framework in which to pose the question `how may $\omega(T) > u$ occur?': we will show that, under $\phniplus$, the optional foundation having length $C \log k$, with $C$ a large constant, is enough to ensure $\omega(T) > u$, and that the conditional probability of this outcome has polynomial decay in $k$. The  $\phniplus$-probability of
$\big\vert V(O_i) \big\vert = k$ will be shown to have exponential decay in $k$. Thus, a logarithmically long optional foundation offers a much more probable alternative to a $k$-sized $i$-th offshoot for realizing the event $\omega(T ) > u$.
\begin{lemma}\label{lemimp}
There exists a constant $c > 0$ such that, for each $\ell \in \N$,
$$
\phniplus \Big( \big\vert V(O_i )\big\vert  = \ell \Big) \leq \exp \big\{ - c \ell \big\}.
$$
There exists a constant $C > 0$
such that
$$
 \phniplus \Big( \omega(T) > u  \Big) \geq k^{-C},
$$
(where recall that $k$ is specified by $\phniplus$). 
The constants $c$ and $C$ may be chosen uniformly in the data that specifies $\phniplus$.
\end{lemma}
We will firstly establish:
\begin{lemma}\label{lemofsc}
Let $i \in\N$. There exists $k_0 \in \N$ such that, for any $i$-absent verdant foliage for which the implied $k$ satisfies $k \geq k_0$, and for any constant $C > 1/(\log \binf)$, the following holds. If a tree $T$ having the given $i$-absent verdant foliage satisfies 
 $\big\vert V(\opf) \big\vert \geq C \log k$, then $\omega(T) > u$.
\end{lemma}
\noindent{\bf Remark on notation.} 
The tree $T$ in the statement of Lemma \ref{lemofsc} has the form of the concatenation of
$\ef$ with $\opf$ and then with $T_{\fbp}$. (By concatenation, we mean the operation under which the endpoint of one path is identified with the root of a second path or tree.) 
The unknown quantities are $O_i$, which determines part of $T_\fbp$, and $\opf$.
We write $T_\fbp = T_\fbp(O_i)$ and $T = T \big( \opf,O_i \big)$ to indicate this dependence. This use of
$T(\cdot,\cdot)$ is in conflict with that of Definition \ref{defnot}, and replaces it from now on.

Note that, for a given tree $O$, the tree $T_\fbp(O)$ may be obtained from $T_\fbp(\emptyset)$
by identifying $\chi_i \in V \big( T_\fbp(\emptyset) \big)$ with the root $\phi_O$ of $O$. For some choices of $O$, the quantity  $T_{\fbp}(O)$ may be ill-defined, however, 
since the choice $O_i = O$ may be incompatible with the given $i$-absent foliage. 
(It may be impossible for  $O$ to play the role of $O_i$, since this may create a vertex whose weight exceeds that of the maximum vertex of $T$.) In the ensuing proof, we find it convenient to work with such trees as $T_\fbp(L_k)$, even though they may technically be illegimate. In such cases, the definition of $T_\fbp(L_k)$ is simply taken to be as in the first sentence of this paragraph.
(Recall that the tree $L_k$ is introduced in Definition \ref{deflkellk}.) \\
\noindent{\bf Proof of Lemma \ref{lemofsc}.} 
% It is in this proof that we invoke the property that $\vmax$ has maximal weight: in brief, the removal of a $k$-sized $O_i$ from $T$ reduces its weight by at most $\frac{k}{k+1}\omega(T)$, since $\vmax$ is not removed; the extension of the foundation by $\big\vert V(\opf) \big\vert \geq C \log k$ is then enough to recover this loss of weight.
We claim that
\begin{equation}\label{claimone}
 \omega^{T(\emptyset,L_k)} \Big( T_\fbp \big( L_k \big) \Big) \geq \frac{\binf - 1}{\binf} u.
\end{equation}
To prove this, note that the expression in the large brackets in (\ref{impineq}) is equal to
$$
\sum_{v \in V \big( T_{\fbp}(L_k) \big)} \omega_{\fbp}(v),
$$
or, equivalently, $\omega^{T_{\fbp}(L_k)}\big( T_\fbp(L_k) \big)$,
the weight of the descendent tree $T_{\fbp}$, rooted at $\fbp$, when the choice $O_i = L_k$ is made.

Thus,
$$
 \omega^{T(\emptyset,L_k)} \Big( T_\fbp \big( L_k \big) \Big) = 
\omega \big( \vnec \big) \sum_{v \in V \big( T_{\fbp}(L_k) \big)} \omega_{\fbp}(v) \geq \frac{\binf - 1}{\binf} u, 
$$
which is (\ref{claimone}).

The equality depends on $\vnec = \fbp$, which is true because $\opf = \emptyset$. The inequality holds because (\ref{impineq}) is verified by $v = \vnec$, by the definition of $\vnec$.

We now claim that
\begin{equation}\label{claimtwo}
  \omega^{T_\fbp(\emptyset)} \Big( T_\fbp \big( \emptyset \big) \Big) \geq \frac{1}{k + 1} 
\omega^{T_\fbp(L_k)} \Big( T_\fbp \big( L_k \big) \Big).
\end{equation}
Indeed, 
$$  
\omega^{T_\fbp(L_k)} \Big( T_\fbp \big( L_k \big) \Big) -
\omega^{T_\fbp(\emptyset)} \Big( T_\fbp \big( \emptyset \big) \Big)
$$
equals $\sum_{v \in V(S_i), v \not= \chi_i} \omega_\fbp(v)$, given that the choice $S_i = L_k$ is made.

Recall that
the offshoot $O_s$ comprises one vertex, which is $\vmax$.
Continuing to set $S_i = L_k$,
we thus have that
\begin{eqnarray}
  & & \sum_{v \in V(S_i), v \not= \chi_i} \omega_\fbp(v) \leq \frac{\vert V(S_i) \vert - 1}{\vert V(S_i) \vert} \Big( \sum_{v \in V(S_i), v \not= \chi_i} \omega_\fbp(v) \,  + \, \omega_\fbp \big( \vmax \big) \Big) \nonumber \\
 & \leq & \frac{\vert V(S_i) \vert - 1}{\vert V(S_i) \vert} 
\omega^{T_\fbp(L_k)} \Big( T_\fbp \big( L_k \big) \Big), \nonumber
\end{eqnarray}
the first inequality since there are $\vert V(S_i) \vert$ summands appearing inside the bracket on its right-hand-side, of which $\omega_\fbp \big( \vmax \big)$ is the largest. It is in deriving this bound that we invoke the defining property of $\vmax$.

We find that
\begin{eqnarray}
& & \omega^{T_\fbp(L_k)} \Big( T_\fbp \big( L_k \big) \Big) -
\omega^{T_\fbp(\emptyset)} \Big( T_\fbp \big( \emptyset \big) \Big) \nonumber \\
& \leq & \frac{k}{k+1} 
\omega^{T_\fbp(L_k)} \Big( T_\fbp \big( L_k \big) \Big), \nonumber
\end{eqnarray}
which yields (\ref{claimtwo}).

We also claim that, for any $m \in \N$, 
\begin{equation}\label{claimthree}
 \omega \Big( T\big( \ell_m,\emptyset \big) \Big) >  \frac{\binf^m}{k+1} \frac{\binf - 1}{\binf} u.
\end{equation}
To this end, note that
$$
\omega \Big( T\big( \ell_m,\emptyset \big) \Big) > \omega^{T ( \ell_m,\emptyset )} \big( \fbp \big)
 \sum_{v \in V\big( T_{\fbp}(\emptyset) \big)} \omega_{\fbp}(v).
$$
Note that, in the case of the tree $T\big( \ell_m,\emptyset \big)$,
$\omega^{T ( \ell_m,\emptyset )} \big( \fbp \big) = \omega \big( \vnec \big) \binf^m$. 
The inequality (\ref{impineq}) holds with $v = \vnec$.
Hence, by (\ref{claimtwo}), we have (\ref{claimthree}):
\begin{eqnarray}
& & \omega \Big( T\big( \ell_m,\emptyset \big) \Big) > \omega^{T ( \ell_m,\emptyset )} \big( \fbp \big)
\omega^{T_\fbp(\emptyset)} \Big( T_\fbp \big( \emptyset \big) \Big) \nonumber \\
 & \geq &  \frac{\binf^m}{k+1} 
\omega \big( \vnec \big)
\omega^{T_\fbp(L_k)} \Big( T_\fbp \big( L_k \big) \Big) 
 \geq \frac{\binf^m}{k+1} \frac{\binf - 1}{\binf}  u, \nonumber
\end{eqnarray}
the third inequality by 
$$
\omega \big( \vnec \big)
\omega^{T_\fbp(L_k)} \Big( T_\fbp \big( L_k \big) \Big)  =
\omega^{T(\emptyset,L_k)} \Big( T_\fbp \big( L_k \big) \Big), 
$$
and by (\ref{claimone}).

By (\ref{claimthree}), $C > 1/(\log \binf)$ and $k \geq k_0$, we have that
$$
\omega \Big( T\big( \ell_{\lfloor C \log k \rfloor},\emptyset \big) \Big) \geq u.
$$
Noting that, for each $m \in \N$, $\omega \big( T(OF,\emptyset) \big)$ is minimized among 
 those $OF$ satisfying $\vert V(OF) \vert = m + 1$ by the choice $OF = \ell_m$,
we  complete the proof of Lemma \ref{lemofsc}. $\Box$

In the preceding lemma, we established a sufficient condition for $\omega(T) > u$
  in terms of the optional foundation $\opf$. For the measures $\phniplus$, we now bound from below the probability that the condition is satisfied,
and find an upper bound on the probability that
$\big\vert V(S_i) \big\vert = k$.

\begin{definition}
Set $q_i = \frac{h_{i+1}}{1 - h_0}$ for $i \geq 0$. Note that $q_i$ is the probability of $i+1$ offspring for   the offspring distribution $\{ h_i: i \in \N \}$ conditioned on there being at least one child. Let $\Theta$ denote the law on weighted trees in which the root has a $\{ q_i: i \in \N \}$-distributed number of offspring, all other vertices having an independent $\{ h_i: i \in \N \}$-distributed number of offspring, and where edge-weights are independently assigned according to the law $\nu$. 
\end{definition}
\begin{lemma}\label{lembdcp}
The distribution of $\opf$ under $\phniplus$ coincides with the distribution of $F$ under $\P_{h,\nu}$.
As such, for $k \geq 1$,
\begin{equation}\label{eqopfb}
\phniplus \Big( \big\vert V\big( \opf \big) \big\vert = k \Big) = h_1^{k-1}\big(1- h_1 \big).
\end{equation}
The distribution of $O_i$ under $\phniplus$ is given by 
a sample tree $R$ from $\Theta$ given that
\begin{equation}\label{weightcond}
 \max_{v \in V(R)} \omega^R(v) < \omega_{\chi_i} \big( \vmax \big).
\end{equation}
There exists $c > 0$ such that, for $i \in \N$, for all $i$-absent verdant foliages, and for $\ell \in \N$, 
\begin{equation}\label{voiexp}
\phniplus \Big( \big\vert V(O_i) \big\vert = \ell  \Big) \leq \exp \big\{ - c\ell  \big\}.
\end{equation}
\end{lemma}
\noindent{\bf Remark.}
Regarding (\ref{weightcond}), note that, under the law $\phniplus$, $\omega_{\chi_i}(\vmax)$ is a deterministic quantity, since it is determined by the biases of the edges in the spine $S$, and these form part of the data of an $i$-absent verdant foliage. \\
\noindent{\bf Proof.} Under the conditional law $\phniplus$, 
the descendent tree $T_{\vnec}$ has the form of the law $\P_{h,\nu}$ conditional on $T_\fbp$
having a collection of offshoots that are consistent with the $i$-absent verdant foliage specified by the choice of $\phniplus$. In the latter law, the length  $\vert E(F) \vert$  of the foundation $F$ does not depend the form of the data in the conditioning, so that this length has the same law as under $\P_{h,\nu}$ conditioned on the sample $T$ being a path. Hence, we obtain (\ref{eqopfb}), since it is clearly true for 
$\P_{h,\nu} \big( \cdot \big\vert \, \textrm{$T$ is a path} \big)$. 
 
To obtain (\ref{weightcond}), note that,
under $\phniplus$, at least one offspring of $\chi_i$ is known to exist, this being $\chi_{i+1}$. Any other offspring form part of $V(O_i)$, by definition. Any form of $O_i$ is admissible, provided that it is compatible with the given $i$-absent verdant foliage. The operative part of this requirement is that no vertex in $O_i$ have a weight exceeding that of the maximum vertex $\vmax$, whose location and weight forms part of the definition of the $i$-absent verdant foliage. The condition that no such vertex exist in $V(O_i)$ is (\ref{weightcond}). 

The assertion (\ref{voiexp}) arises as follows. 
It follows directly from Lemma \ref{sizedec} that
$$
\Theta \Big( \big\vert V(R) \big\vert = \ell \Big) \leq \exp \big\{ - c \ell \big\}
$$
for all $\ell \in \N$, 
(where $R$ denotes the sample of the measure $\Theta$).
Thus, (\ref{voiexp}) follows from
\begin{equation}\label{eqwc} 
\Theta
 \Big(
 \max_{v \in V(R)} \omega^R(v) < \omega_{\chi_i} \big( \vmax \big)
  \Big) \geq c.
\end{equation}
We now explain why (\ref{eqwc}) holds. 
Writing $\{ \phi \}$ for the singleton graph consisting only of a root, we see that the choice 
$R = \{ \phi \}$  implies (\ref{weightcond}), since it reduces to the statement $1 < \omega_{\chi_i} \big( \vmax \big)$, which is true because $\chi_i \not= \vmax$ is implied by the definition of an $i$-absent foliage (in the form of the assumption that $i < s$). However, $h_1 > 0$ implies that $\Theta \big( R = \{ \phi \} \big) > 0$, whence (\ref{eqwc}). \\
%\begin{lemma}
%$$
%\P_{h,\nu}^+ \Big( \omega(T) > u \Big) \geq k^{-C}
%$$
%and
%$$
%\P_{h,\nu}^+ \Big( \big\vert V(O_i) \big\vert = k  \Big) \leq \exp \big\{ - ck  \big\}.
%$$
%\end{lemma}
\noindent{\bf Proof of Proposition \ref{lemsmallos}.}
We have that
\begin{eqnarray}
& & \phni \Big(  \omega(T) > u, \big\vert V(O_i) \big\vert = k  \Big) 
 \nonumber \\
 & = & \phni \Big(  \ess  \Big)
\phni \Big(  \omega(T) > u, \big\vert V(O_i) \big\vert = k  \Big\vert \ess  \Big) \nonumber \\
 & \leq & \phni \Big(  \ess  \Big) \exp \big\{ - ck \big\}.
 \nonumber
\end{eqnarray}
The inequality here follows by conditioning on the form of $EF$ and then applying the final assertion of Lemma \ref{lembdcp}.
On the other hand,
\begin{eqnarray} 
& & \phni \Big(  \omega(T) > u \Big) 
 \geq  
 \phni \Big(  \ess  \Big)
\phni \Big(  \omega(T) > u  \Big\vert \ess  \Big) \nonumber \\
 & \geq & \phni \Big(  \ess  \Big) \phniplus \Big( \big\vert V\big(\opf\big) \big\vert \geq C \log k \Big) 
  \geq   
\phni \Big(  \ess  \Big) k^{- 2 C \log \big( h_1^{-1} \big)}. \nonumber 
\end{eqnarray}
In the use of the law $\phniplus$, an arbitrary form for the fixed essential foundation $\ef'$ may be taken.
The inequalities are due to Lemma \ref{lemofsc} and to Lemma \ref{lembdcp}(i).
By the two preceding displayed statements,
$$
\phni \Big(  \big\vert V(O_i) \big\vert = k \Big\vert \omega(T) > u \Big) \leq k^{ 2C \log \big( h_1^{-1} \big)} \exp \big\{ - ck \big\}.
$$ 
By decreasing the value of $c > 0$, we obtain the statement of the lemma (under Hypothesis \ref{hyplem}).

It remains to explain how the argument for Proposition \ref{lemsmallos}
may be modified so that Hypothesis \ref{hyplem} is not invoked. 

Firstly, suppose that the offspring law satisfies $h_1 = 0$. 
Set $\kappa = \inf \big\{ \ell \geq 2: h_\ell > 0 \big\}$. (Note that $\kappa < \infty$, since we made have made this assumption from the outset to avoid trivialities.) With its present definition, the foundation $F$ is necessarily empty.
We alter the definition as follows. In a rooted tree $T$, define a vertex $v$ to be standard if the path $P_{\phi,v}$ 
has the property that $u \in V\big(P_{\phi,v}\big)$, $u \not= v$, implies that $u$ has $\kappa$ offspring;
that all of the $\kappa - 1$ of these offspring not lying in $V\big(P_{\phi,v}\big)$ have no offspring;
and that, of the weights attached to the $\kappa$ edges connecting $u$ to its offspring, there is a unique largest, and
it is associated to the edge connecting $u$ to its offspring in $V\big(P_{\phi,v}\big)$. Note that $\phi$ is standard.

We now define $\fbp$ to be the final vertex in the initial sequence of standard vertices in the path 
$P_{\phi,\vmax}$. The foundation is then taken to be the 
connected component of $\phi$ in the graph with vertex set $V(T)$ and edge-set $E(T) \setminus G_\fbp$,
where $G_\fbp$ denotes the set of edges connecting $\fbp$ to its offspring (that is, the foundation is taken to be the part of $T$ ``above'' $\fbp$). 

Certain changes in the argument are forced. We comment only on the principal ones. The weighted trees $L_k$ and $\ell_k$ in Definition \ref{deflkellk} are now defined to be of depth $k$, with vertices 
at distance from the root less than $k$
each having $\kappa$ offspring. 
The requirement on edge-biases in the new definition of foundation 
ensures that an arbitrary choice of optional foundation $\opf$ is compatible with a given $i$-absent verdant foliage,
since it implies that vertices $v$ lying in $\opf$ but not in $P_{\phi,\vmax}$ have 
$\omega(v) < \omega \big(\vmax  \big)$. In deriving the final assertion (\ref{voiexp}) of Lemma \ref{lembdcp},
we must show (\ref{eqwc}). 
We may no longer consider the event that $R = \{ \phi \}$, but instead consider the event 
$F_{\kappa-1}$ 
that $R$ consists of a root with  $\kappa - 1$ offspring, with these vertices each having no offspring. This event has $\Theta$-probability at least $c > 0$. To show (\ref{eqwc}), we need to confirm that
\begin{equation}\label{evinc}
F_{\kappa - 1} \subseteq \Big\{ \omega^R(v) < \omega_{\chi_i}(\vmax) \, \, \textrm{for all $v \in V(R)$} \Big\}.
\end{equation}
To do so, we note that, if there is no tree $R$ of depth $\ell$ 
such that  $\omega^R(v) < \omega_{\chi_i}(\vmax)$ for all $v \in V(R)$, 
then (\ref{voiexp}) certainly holds, since its left-hand-side is zero. 
This being what we seek to show, we may assume the other case, and let $T'_i$ denote some such tree. Provided that $\ell \geq \frac{\log \bsup}{\log \binf}$, any choice of $R$ in the support of $F_{\kappa - 1}$ has the property that each element of $V(R)$ has weight less than that of some vertex in $T'$. This confirms (\ref{evinc}), so that we have (\ref{voiexp}) for all $\ell$ sufficiently high. We obtain the general statement by decreasing the value of $c > 0$.

We also suspend the assumption that $\nu$ have no atoms. In defining the $\fso$-decomposition, we must decide which of the vertices of maximal weight $\vmax$ should be. For ordered trees, we might choose $\vmax$ to be the lexicographically minimal vertex among those of maximal weight. However, as we discussed after the statement of Proposition \ref{lemsmallos}, we may not work with ordered trees.
 Our solution is to use the lexicographical ordering of a sample of $\P_{h,\nu}$
merely to identify $\vmax$, after which, we forget about this ordering. This is the formal definition:
\begin{definition}
For 
a finite rooted ordered weighted tree $T$,
we write $\vmax$ for the lexicographically minimal element of $V(T)$ among those of maximal weight.
Two such trees $T$ and $T'$ are said to be $\vmax$-isomorphic
if there exists an isomorphism $\psi:V(T) \to V(T')$ in the sense of Definition \ref{defformphnu} with the property that
$\psi(\vmax^T) = \vmax^{T'}$.

A finite rooted unordered weighted tree {\it with declared maximum vertex} is an equivalence class of such trees under the equivalence relation of $\vmax$-isomorphism. 

An unordered weighted Galton-Watson tree with {\it declared maximum vertex}
is the law on finite rooted unordered weighted trees with declared maximum vertex
given by the $\vmax$-isomorphism class of a sample of an ordered tree under $\P_{h,\nu}$.
%arising from choosing $\vpoint = \vmax$, where $\vmax$ is the lexicographically minimal vertex of maximal weight in $V(T)$. 
% is said to be {\it pointed} if associated to $T$ is a distinguished element $\vpoint \in V(T)$.  
%A finite rooted ordered edge-weighted tree $T$ is said to be {\it pointed} if associated to $T$ is a distinguished element $\vpoint \in V(T)$. 
%Two such trees $T,T'$ are said to be pointed isomorphic if there exists an isomorphism $\psi:V(T) \to V(T')$ in the sense of Definition \ref{} with the firther property that $\psi\big( \vpoint_T \big) = \vpoint_{T'}$.
%A finite rooted unordered edge-weighted pointed tree is an pointed isomorphism class of such trees.
%An unordered edge-weighted Galton-Watson tree with {\it declared maximum vertex}
%is the law on finite rooted unordered edge-weighted pointed trees
%arising from choosing $\vpoint = \vmax$, where $\vmax$ is the lexicographically minimal vertex of maximal weight in $V(T)$. 
\end{definition}
We now suspend the notational abuse by which $\P_{h,\nu}$ denoted the unordered Galton-Watson law. For the ensuing argument, we will abuse notation in a different way, writing $\P_{h,\nu}$ for the unordered weighted Galton-Watson tree with declared maximum vertex.
It is for this new choice of $\P_{h,\nu}$ that the lemmas in the preceding proof are to be understood, in the case that $\nu$ has no atoms.

%In considering the altered definition of an $i$-absent (verdant) foliage, it would appear natural to include in the definition of the spine and the off-shoots, a natural number associated to each edge, indicating the order of the edge among the offspring of its parent vertex. However, this results in too much conditioning, since, if the edge $[\chi_i,\chi_{i+1}]$ has a high index, this forces $O_i$ to have a root with many children, so that $O_i$ will necessarily be large.

%To circumvent this difficulty, we do not include this ordering data in the definition of  an $i$-absent (verdant) foliage. 

Regarding the statement of Lemma \ref{lembdcp}, let $S$ denote the distribution of $O_i$ under $\phniplus$. We now have that, for $R$ an unordered weighted tree in the support of $\Theta$, the Radon-Nikodym derivative $\frac{dS}{d \Theta}$ satisfies
\begin{equation}
      Z^{-1} \frac{dS}{d \Theta} \big( R \big)  =  
   \left\{ \begin{array}{rl}            
1 & \textrm{if }   
\max_{v \in V(R)} \omega^R(v) < \omega_{\chi_i} \big( \vmax \big),  \\
    p_R  & \textrm{if }   
\max_{v \in V(R)} \omega^R(v) = \omega_{\chi_i} \big( \vmax \big), \\
   0 &  \textrm{otherwise,}  \end{array} \right.
\end{equation}
where $Z > 0$ is a normalization. The quantity $p_R$ has the following interpretation.
Consider the (unordered weighted) tree (with declared maximum vertex) arising from $\phniplus$  with the choice $O_i = R$ being made. Recall that this tree has no ordering, but that it does carry a declared vertex of maximal weight.
Now assign a random ordering to the set of offspring of each vertex in this tree, with the uniform law. 
Then $p_R$ is equal to the probability that the lexicographically minimal vertex of maximal weight, as selected by this random ordering, is equal to the declared vertex. (Note that, in the case that $\max_{v \in V(R)} \omega^R(v) = \omega_{\chi_i} \big( \vmax \big)$, we have that $p_R < 1$, since the copy of $R$ playing the role of $O_i$ contains at least one vertex of maximal weight.)  

In fact, the actual value that $p_R$ takes is irrelevant: our purpose is served by verifying (\ref{voiexp}), and, 
given that (\ref{voiexp}) is valid for $\Theta$, (as follows immediately from Lemma \ref{sizedec}), it remains only to show that $Z \geq c > 0$. However, 
$$
Z \geq \Theta \Big( \max_{v \in V(R)} \omega^R(v) < \omega_{\fbp} \big( \vmax \big) \Big),
$$ 
The right-hand-side was shown to be positive by the argument involving (\ref{evinc}).Hence, $Z > c > 0$, as required for (\ref{voiexp}) in the case that neither element of Hypothesis \ref{hyplem} is assumed.  
This completes the proof of Proposition \ref{lemsmallos}. $\Box$
\end{subsection}
\begin{subsection}{Deriving Theorem \ref{lemdivexpold}  and Lemma \ref{lemaux} from Proposition \ref{lemsmallos}}
\noindent{\bf Proof of Lemma \ref{lemaux}.} Recalling that $s = \vert V(S) \vert - 1$, we label the offshoots $\big\{ O_i: 0 \leq i \leq s \big\}$ appearing in the $\fso$-decomposition of a tree $T$ in reverse order as
$\hat{O}_i = O_{s - i}$, $0 \leq i \leq s$, setting 
$\hat{O}_i = \emptyset$ if $i > s$.

By averaging Proposition \ref{lemsmallos} over the distribution of the $(\vert V(S) \vert - i)$-th absent foliage under $\P_{h,\nu,u}$, we obtain the following consequence.
\begin{lemma}\label{lemconseq}
There exists $c > 0$ such that, for all $i,k \in \N$ and for all $u > 0$,
\begin{equation}\label{eqhato}
  \P_{h,\nu,u} \Big(  \big\vert V\big(\hat{O}_i\big) \big\vert \geq k  \Big) \leq \exp \big\{ - c k \big\}.
\end{equation}
\end{lemma}
Lemma \ref{lemconseq} implies that, for some $c > 0$,
\begin{equation}\label{oiprimebd}
  \sup_{u > 0}
  \P_{h,\nu,u} \Big( \sum_{i=m}^{\vert V(S) \vert} \big(2\bsup\big)^{\big\vert V \big( \hat{O}_i \big) \big\vert}  
  \binf^{-i}  \geq \exp \big\{ - cm \big\} \Big) \leq \exp \big\{ - cm \big\}  
\end{equation}
for each $m \geq 1$.
Let
$$
A_\ell = \Big\{   d \big( \vbase, \vmax \big) \leq \ell/2 \Big\}.
$$
We observe that, on $A_\ell$, the sequence 
$\big\{ J_{D(T) - \ell - i}: i \geq 0 \big\}$
coincides with
$\big\{ \hat{O}_{m+i}: i \geq 0  \big\}$ for some $m \in \big\{ \ell/2,\ldots,\ell \big\}$. We now verify this, with the aid of Figure 3.
%\ref{fighatoj}.
\begin{figure}\label{fighatoj}
\begin{center}
\includegraphics[width=0.6\textwidth]{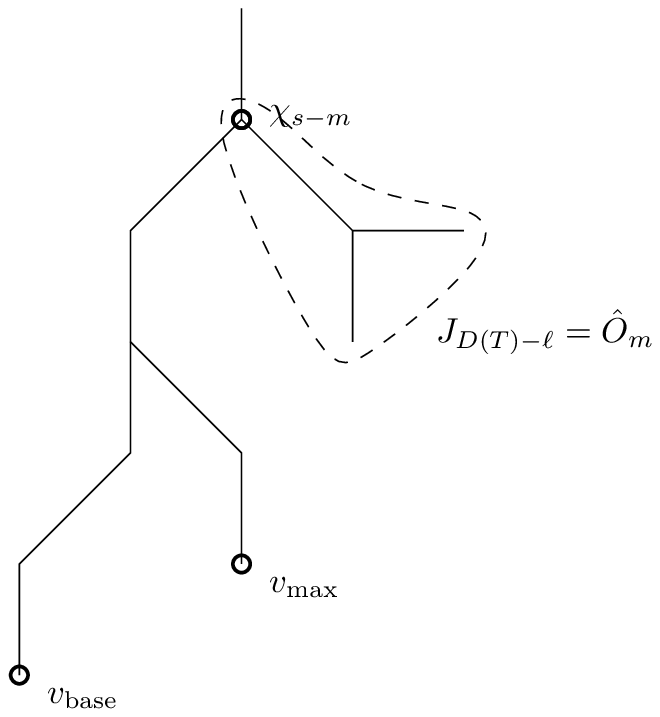} \\
\end{center}
\caption{A tree viewed in the neighbourhood of $\vbase$ and $\vmax$.}
\end{figure} 
The outgrowth $J_{D(T) - \ell}$, whose root is the vertex $v$ lying at distance $\ell$ from $\vbase$, is also the offshoot $\hat{O}_m$, where $m = d \big(\vmax,v \big)$. 
We have that $m \leq \ell$,  
since $d \big( \phi , \vbase \big)  \geq d \big( \phi, \vmax \big)$, while 
 $d \big( \vbase, \vmax \big) \leq \ell/2$ forces $m \geq \ell/2$.

Hence, on $A_\ell$,
$$
 \sum_{i=0}^{D(T) - \ell} \big(2\bsup\big)^{\big\vert V \big( J_{D(T) - \ell - i} \big) \big\vert}
 \binf^{-\ell-i} \leq
  \sum_{i=0}^{\vert V(S) \vert - \ell/2} \big(2\bsup\big)^{\big\vert V \big( \hat{O}_{i + \ell/2} \big) \big\vert}
 \binf^{-\ell/2 - i}, 
$$
again by 
$d \big( \phi , \vbase \big)  \geq d \big( \phi, \vmax \big)$.
By (\ref{oiprimebd}), then,
$$
  \sup_{u > 0}
  \P_{h,\nu,u} \Big( A_\ell, \sum_{i=0}^{D(T) - \ell} \big(2\bsup\big)^{\big\vert V \big( J_{D(T) - \ell - i} \big) \big\vert}  
   \binf^{-i-\ell}  \geq \exp \big\{ - c\ell/2 \big\} \Big) \leq \exp \big\{ - c\ell/2 \big\}.  
$$
For the statement of Lemma \ref{lemaux}, it suffices then to show that
\begin{equation}\label{akcbd}
   \sup_{u>0} \P_{h,\nu,u} \big( A_\ell^c \big) \leq \exp \big\{ - c \ell \big\}.
\end{equation}
Let $\ancestor$ denote the last common ancestor of $\vbase$ and $\vmax$.
We argue that
\begin{equation}\label{maxbase}
 A_\ell^c \subseteq \bigg\{ \Big\vert  V \big( \hat{O}_{d(\vbase,\vmax)} \big) \Big\vert \geq \max \Big\{ d \big( \ancestor, \vmax \big), \ell/4 \Big\} \bigg\}. 
\end{equation}
To see this, 
note that $\vbase \in V\big( \hat{O}_{d(\ancestor,\vmax)} \big)$, whence $\big\vert  V\big( \hat{O}_{d(\ancestor,\vmax)} \big\vert \geq d \big( \ancestor,\vbase \big)$.  
Now, $d \big( \ancestor,\vbase \big)$ is at least $d(\ancestor,\vmax)$, since $d \big( \phi, \vbase \big) \geq d \big( \phi , \vmax \big)$, while it exceeds 
$\ell/4$ on $A_\ell^c$, because then 
$d(\vmax,\ancestor) + d(\ancestor,\vbase) = d(\vmax,\vbase) > \ell/2$.  

From (\ref{maxbase}), we learn that
$$
\P_{h,\nu,u} \big( A_\ell^c \big) \leq 
 \sum_{i = 1}^{\ell/4} \P_{h,\nu,u} \Big( \big\vert V\big(\hat{O}_i \big) \big\vert \geq \ell/4  \Big)
 \, + \,
\sum_{i = \ell/4 + 1}^\infty \P_{h,\nu,u} \Big( \big\vert V\big(\hat{O}_i \big) \big\vert \geq i  \Big).
$$
Applying Lemma \ref{lemconseq} confirms (\ref{akcbd}) and completes the proof of Lemma \ref{lemaux}. $\Box$ \\
We now provide two lemmas in preparation for the proof of Theorem \ref{lemdivexpold}.
As we have noted, Theorem \ref{lemdivexpold} is trivial if the hypothesis of the next result is violated.
\begin{lemma}\label{lemfoundbd}
Assume that $h_i > 0$ for some $i \geq 2$.
The foundation $F$ in the FSO-decomposition satisfies the following bound. 
There exists $c > 0$ such that, for all $k \in \N$,
$$
\sup_{u > 0} \P_{h,\nu,u} \Big( \big\vert V(F) \big\vert \geq k \Big) \leq \exp \big\{ - ck \big\}.
$$ 
\end{lemma}
\noindent{\bf Proof.} 
We consider the case that Hypothesis \ref{hyplem} holds,
the other case requiring only minor modifications.

Consider a procedure by which the law $\P_{h,\nu}$ is constructed by firstly realizing the offspring of the root, and the associated biases on edges, and then iteratively selecting an as-yet-selected vertex and similarly realizing its offspring and the intervening biases (so that the new offspring are added to the list of vertices awaiting selection). The procedure stops when all vertices in the presently constructed tree have been selected; naturally, this stopping occurs in a finite number of steps for a subcritical law $h$. The means of choosing the next vertex to be selected may be any fixed previsible procedure.
Let $Q \subseteq V(T)$ be such that $q \in Q$ if and only if, at the moment the procedure selects $q$, $d(\phi,q)$ is strictly greater than $d(\phi,r)$ for any already selected vertex $r$. Clearly, $\vert Q \vert$ is one plus the depth $D(T)$ of the sampled tree $T$.
Note that, if the event $\big\vert V(F) \big\vert \geq k$ is to occur, it must be the case that, for each of the first $k-1$ elements of $Q$ selected by the procedure, 
the number of offspring realized is equal to one. It is easy to see, however, that, for any $u > 0$ and $i \in \N$, under the conditional law 
$\P_{h,nu,u} \big( \cdot \big\vert \vert Q \vert \geq i \big)$, the number of offspring realized for the $i$-th selected element of $Q$ stochastically dominates the unconditioned offspring distribution $h$. 
This means that, under $\P_{h,\nu,u}$, each successive element of $Q$ may have at least two offspring, with the probability of this being bounded below, conditionally on the construction made thus far; the probability that this fails to happen on $k-1$ consecutive occasions decays exponentially in $k$, whence the result. $\Box$

\begin{lemma}\label{lemopr}
Write $P_{\phi,\vmax} = \big[ \phi = \chi'_0,\chi'_1,\ldots,\chi'_{d(\phi,\vmax)} = \vmax \big]$.
Let $O'_i$ denote the connected component containing $\chi'_i$ in the graph with vertex set $V(T)$
and edge-set $E(T) \setminus E \big( P_{\phi,\vmax} \big)$. Then there exists $c > 0$ such that, for all $i,k \in \N$,
$$
\sup_{u > 0} \P_{h,\nu,u} \Big( \big\vert V (O'_i) \big\vert \geq k \Big) \leq \exp \big\{ - ck \big\}.
$$  
\end{lemma}
\noindent{\bf Proof.}
Note that, in the event that $\big\vert V(F) \big\vert = i + 1$, then $O'_j = O_{j-i}$ for $j \geq i$. Hence,
$$
\Big\{   \big\vert V\big(O'_i\big) \big\vert \geq k \Big\}
\subseteq
 \bigcup_{j = i - k/2}^i 
\Big\{   \big\vert V\big(O_j\big) \big\vert \geq k \Big\} \, \cup \, 
\Big\{   \big\vert V\big( F \big) \big\vert \geq k/2 \Big\}.
$$
Thus,
$$
\P_{h,\nu,u} \Big( \big\vert V\big(O'_i\big) \big\vert \geq k  \Big)
\leq \frac{k}{2}  \sup_j 
\P_{h,\nu,u} \Big( \big\vert V\big(O_j\big) \big\vert \geq k  \Big) \, + \,
\P_{h,\nu,u} \Big( \big\vert V\big(F\big) \big\vert \geq k/2  \Big). 
$$
By Proposition \ref{lemsmallos} and Lemma \ref{lemfoundbd}, we obtain the statement of the lemma.\\
\noindent{\bf Proof of Theorem \ref{lemdivexpold}.}
Recall that $\ancestor$ denotes the last common ancestor of $\vmax$ and $\vbase$.
Note that $\ancestor \in V(P_{\phi,\vmax})$, so that we may write $\ancestor =  \chi'_M$, with $0 \leq M \leq d \big( \phi , \vmax \big)$.

Recall that $i \in \N$ denotes the index fixed in the statement of the theorem. 
Set $\mathcal{M}_i = \big\{M > i \big\}$.

Note that $J_j = O'_j$
 for $j \leq M - 1$, 
so that $\mathcal{M}_i \subseteq \big\{  J_i = O'_i \big\}$. Note further that 
$\mathcal{M}_i^c \subseteq \big\{  V(J_i) \subseteq  V(T_{\ancestor}) \big\}$.

We have that
\begin{eqnarray}
& & \P_{h,\nu,u} \Big(  \big\vert  V \big( J_i \big) \big\vert \geq k \Big) \nonumber \\
 & \leq & 
\P_{h,\nu,u} \Big(  \big\vert  V \big( J_i \big) \big\vert \geq k , \mathcal{M}_i \Big)
 + 
\P_{h,\nu,u} \Big(  \big\vert  V \big( J_i \big) \big\vert \geq k , \mathcal{M}_i^c \Big) \nonumber \\
  & \leq &  
\P_{h,\nu,u} \Big(  \big\vert  V \big( O'_i \big) \big\vert \geq k  \Big)
 + 
\P_{h,\nu,u} \Big(  \big\vert  V \big( T_\ancestor \big) \big\vert \geq k  \Big).
\end{eqnarray}
Hence, the statement of the theorem follows from Lemma \ref{lemopr} and the next assertion. 

There exists $c > 0$ such that, for all $u > 0$, and $k \in \N$,
\begin{equation}\label{eqvtanc}
  \P_{h,\nu,u} \Big(  \big\vert  V \big( T_\ancestor \big) \big\vert \geq k  \Big) \leq \exp \big\{ - ck \big\}.
\end{equation}
To prove (\ref{eqvtanc}), note that
$$
V \big( T_\ancestor \big) \subseteq \bigcup_{j=0}^{d \big( \vmax,\vbase \big)} V \big( \hat{O}_j \big).
$$
Hence, (\ref{eqvtanc}) is a consequence of Lemma \ref{lemconseq} and (\ref{akcbd}). $\Box$
\end{subsection}
\end{section}
\begin{section}{The renewal decomposition of a rooted tree}\label{appc}
The tree decompositions that we have inroduced have been sufficient for the purposes of proving the results in this paper. 
In this appendix, we introduce a further decomposition of a rooted tree, whose component enjoys an attractive independence property. 
This ``renewal decomposition'' splits the tree at a sequence of points that have a natural interpretation as regeneration points. This structure permits us to describe explicitly the conditional distribution of $\P_{h,\nu}\big( \cdot \big\vert \omega(T) > u \big)$, given the form of a finite initial sequence of components from the renewal decomposition. In Proposition \ref{lemdivexp}, we will prove the analogue of Theorem \ref{lemdivexpold} for the new decomposition.
In its derivation of a stable limiting law for randomly biased walk on a supercritical Galton-Watson tree,
the paper \cite{GerardAlan} will make an essential use of the renewal decomposition, and of Proposition \ref{lemdivexp}. 
\begin{subsection}{Definition of the renewal decomposition}
\begin{definition}\label{defndiv}
By a root-base tree $T$, we refer to a finite rooted tree, one of whose vertices $\basedef$ at maximal distance from $\phi$ is declared to be the base.

Given a rooted tree $T$, a vertex $v \in V(T)$, $v \not= \phi$, is called a cutpoint if   
it is not a leaf, and any other vertex in $T$ at the same distance from $\phi$ as $v$ is a leaf.
The set of cutpoints naturally decompose a rooted tree into components in the following manner. We write $r(T)$ for the number of cutpoints of $T$ plus one.
We may then record these cutpoints in the form $c_i$, $1 \leq i \leq r(T) - 1$, in increasing order of distance from the root $\phi$. 
We further set $c_0 = \phi$. 
We set $d_i = d \big( \phi, c_i \big)$ for 
$0 \leq i \leq r(T) - 1$. We also set $d_{r(T)}= D(T)$, where recall that $D(T) = \max \big\{ d \big( \phi,v \big): v \in V(T) \big\}$.  
For $1 \leq i \leq r(T)$, we define the $i$-th component $C_i$ of the tree $T$ to be the subgraph of $T$ induced by the set of vertices in $T$ at a distance from the root of at least $d_{i-1}$ and at most $d_i$.
Then, for $1 \leq i \leq r(T) - 1$, $C_i$ may be regarded as a root-base tree, with
$\phi(C_i) = c_{i-1}$ and $\basedef(C_i) = c_i$. The final component $C_{r(T)}$, however, is a tree rooted at $c_{r(T) - 1}$ that has no natural choice of base. 
\end{definition}
In this way, we divide a rooted tree into components. We may also perform the operation in reverse, assembling a collection of such trees into a single one.
\begin{definition}\label{defcirc}
Let $T_1,\ldots,T_r$
be a finite sequence of finite rooted trees, 
all but the last of which is root-base.
We define the concatenation
$$
T_1 \circ \cdots \circ T_r
$$
to be the rooted tree with vertex set $\bigcup_{i=1}^r V(T_i)$,
in which the identifications $\basedef(T_i) = \phi(T_{i+1})$, $1 \leq i \leq r-1$,
are made, and whose edges are inherited from the constituent trees.
The root of $T_1 \circ \cdots \circ T_r$ is taken to be equal to $\phi(T_1)$.
\end{definition}
\noindent{\bf Remark.}
Note that the components $\{ C_1,\ldots,
C_{r(T)}\}$ of a finite rooted
tree $T$ appear as the constituents in the decomposition
$T = C_1 \circ \ldots \circ C_{r(T)}$, in which identifications between
successive elements are made at the cutpoints of $T$. 
\end{subsection}
\begin{subsection}{Sampling the components of a subcritical tree}
To gain an understanding of the renewal decomposition, 
we will form a sample of $\P_h$ iteratively, adding one component
at a time. The following definition is convenient.
\begin{definition}
Let $\P_h^*$ denote the law of $\P_h$ given that the sample $T$ is non-trivial, that is, given that $E(T)$ contains at least one element. We further write $\P_{h,\nu}^*$ for the law on rooted weighted trees in which the edges of a sample of $\P_h^*$ are independently assigned biases, each having law $\nu$.
\end{definition} 
Having constructed the first $k$ components of a sample of $\P_h^*$, we will add the
$(k+1)$-st according to the correct conditional distribution. A subtlety
arises, because the final component lacks a base. 
%In speaking of the law of
%descendent tree $T_{b(C_k)}$ under $\P_h^*$, given
%the values of its first $k$ components $C_1, \ldots, C_k$,
%then, we must bear in mind that we implicitly presuppose that $r(T) \geq k+1$, in order that the base %$b(C_k)$
%may be well-defined. The precise statement is that, 
Note that, for each $k \in \N$, 
under the law $\P_h$ conditional on $r(T) \geq k+1$ and the values of
$C_1,\ldots,C_k$, the tree $T_{\basedef(C_k)}$ has law $\P_h^*$.
The next component $C_{k+1}$ is root-base precisely when it is not the final
component, namely, with probability  $\P_h^* \big( r(T) \geq 2 \big)$, which is the conditional probability that a non-trivial tree contains a cutpoint; in this event, the distribution of $C_{k+1}$ is that of $C_1$
under  $\P_h^* \big( \cdot \big\vert r(T) \geq 2 \big)$, and it has a base to which a further component
will be added; otherwise, it has the law of $C_1$ under $\P_h^* \big( \cdot
\big\vert r(T) = 1 \big)$. 
We now present two lemmas that summarise these conclusions.
\begin{lemma}\label{lemdecom}
Let $\rho_{\rm r,b}$ denote the law on root-base trees given
by $C_1$ under the law $\P_h$ conditional on $T$ containing a cutpoint. 
Let $\rho_{\rm r}$ denote the law on rooted 
trees given by $C_1$ under  $\P_h$ conditional on $T$ being non-trivial but containing no cutpoint. 

Then the decomposition $T = C_1 \circ \ldots \circ C_{r(T)}$
 under $\P_h^*$ has the form 
$$
T = C_1 \circ \ldots \circ C_{L+1},
$$
with $L$ an $\N$-valued random variable satisfying
$$
\P \big( L = k \big) = \P_h^* \big( r(T) \geq 2 
\big)^k   \P_h^* \big( r(T) = 1 \big), \qquad k \geq 0,
$$
and with $C_1,\ldots C_L$ denoting a sequence of independent root-base
trees, each having law $\rho_{\rm r,b}$, while $C_{L+1}$ is an independent
rooted tree with law $\rho_{\rm r}$. $\Box$
\end{lemma}
\begin{lemma}\label{lemdecomnew}
Let $k \in \N$.
Consider the law $\P_{h,\nu}$ conditionally on $r(T) \geq k+1$ and on an arbitrary form for the first $k$ weighted components of $T$. Then the conditional distribution of the descendent tree $T_{\basedef(C_k)}$
is given by $\P_{h,\nu}^*$.
\end{lemma}
Contrast Lemma \ref{lemdecomnew} with the analogous result for the outgrowth splitting of a rooted tree introduced in Definition \ref{defoutgr}. If we condition on $D(T) \geq k + 1$, 
and on the form of the first $k+1$ outgrowths $\big\{J_i: 0 \leq i \leq k \big\}$, then the conditional distribution of $T_{\psi_{k+1}}$
is dependent on the form of these initial outgrowths, because of the requirement that the vertex $\vbase$ lie in
$T_{\psi_{k+1}}$. (The vertex $\vbase$ being of maximal depth, it is necessary that $T_{\psi_{k+1}}$ be deep enough to reach further from the root of $T$ than all of the initial conditioned outgrowths.)
As Lemma \ref{lemdecomnew} shows, this difficulty does not arise for the renewal decomposition. For this reason, it is natural to think of the cutpoints of in the renewal decomposition
as ``regeneration'' points for the tree. The renewal decomposition, in the form of Lemma \ref{lemdecomnew}, will be exploited in \cite{GerardAlan}, to understand the geometry of a subcritical trap, conditional on its structure near the trap entrance.

Lemma \ref{lemdecom} has the following corollary. Recall that we use the notation
 $P_{h,\nu,u} = P_{h,\nu} \big( \cdot \big\vert \omega(T) > u \big)$.
\begin{corollary}\label{cordecom}
Let $k \in \N$ and $u > 0$.
Consider the law $\P_{h,\nu,u}$ conditionally on $r(T) \geq k+1$ and on an arbitrary form for the first $k$ weighted components of $T$.  Set $u' \in \R$ according to
$$
 u' = \frac{u - \sum \Big\{ \omega(v) : v \in \bigcup_{i=1}^k V \big( C_i \big) \setminus \big\{ \basedef(C_k) \big\} \Big\}
 }{\omega \big( \basedef(C_k) \big)}.
$$
If $u' > 0$, then the conditional distribution of
$T_{\basedef(C_k)}$ has the
law $\P_{h,\nu}\big( \cdot \big\vert \omega(T) > u'  \big)$;
if $u' \leq 0$, it has the law $\P_{h,\nu}^*$.
\end{corollary} 
%The explicit form of this conditional distribution, uncluttered by additional biases found in the case of other %tree-splittings, 
%indicates how the cutpoints of a rooted tree
%may be considered as ``regeneration'' points for the tree.  
\end{subsection}
\begin{subsection}{Components in a high-weight tree are small}
We now present the analogue of Theorem \ref{lemdivexpold} for the renewal decomposition.
\begin{prop}\label{lemdivexp}
For the statement, we take $C_i = \emptyset$ if $i > r(T)$ (for any rooted
tree $T$). There exists $c > 0$ such that, for all $u > 0$ and $i \in \mathbb{N}$,
$$
\mathbb{P}_{h,\nu,u} \Big( \big\vert V \big( C_i \big) \big\vert \geq k \Big) \leq \exp \big\{ - c k  \big\},
$$
for each $k \in \mathbb{N}$.
\end{prop}

Our task is to show that the component sizes of a rooted weighted tree have exponential tails, uniformly under conditioning
a sample $T$ of $\P_{h,\nu}$ on $\omega(T) > u$.
Firstly, consider Figure 4,
%\ref{figtwosketches}, 
in which two problematic trees are depicted. 
Recalling Definition \ref{defdecomplem}, in each sketch, the vertical segment represents the path $P_{\phi,\vmax}$, and the slanting segments the offshoots. Each tree has a large component (in its renewal decomposition), indicated in a dashed box. In the left-hand sketch, the cause of the large component is a long outgrowth, so that the scenario is excluded in essence by Proposition \ref{lemsmallos}.  In the right-hand sketch, the cause occurs is a succession of overlapping short outgrowths.
\begin{figure}\label{figtwosketches}
\begin{center}
\includegraphics[width=0.7\textwidth]{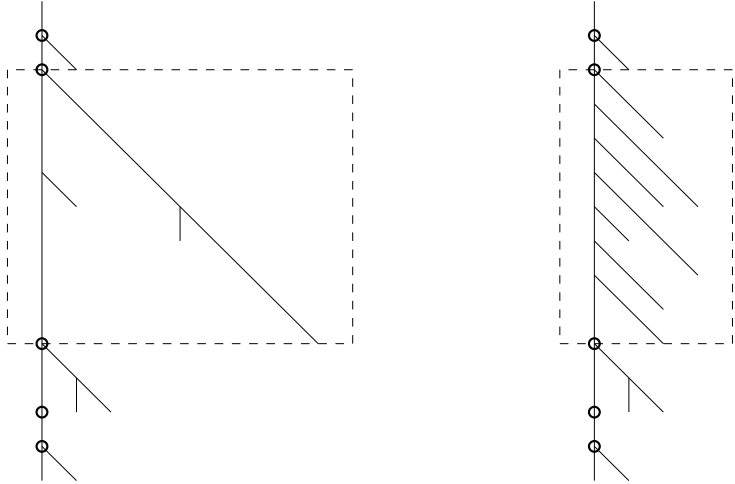} \\
\end{center}
\caption{Two trees having a large component. The small circles indicate cutpoints.}
\end{figure}
Not only does Proposition \ref{lemsmallos} establish that a given offshoot is unlikely to be long, it does so conditionally on the lengths of the other offshoots. Note the dependence in the offshoot lengths in the second sketch in Figure 4.
\ref{figtwosketches}. 
Proposition \ref{lemsmallos} will also be valuable in showing that the right-hand  scenario is improbable, since it entails much dependence in the lengths of successive offshoots.

Alongside Proposition \ref{lemsmallos}, 
The following result will be required.
\begin{lemma}\label{lememptyos}
There exists $c > 0$ such that, for any $u > 0$,
and  for all $i \in \N$,
$$
\phni \Big(  D(O_i)  \leq 1 \Big\vert  \omega(T)> u \Big) \geq c.
$$ 
\end{lemma}
\noindent{\bf Proof.} Consider firstly the case that $h_1  > 0$. The details of this argument are a reworking of those for Proposition \ref{lemsmallos}, and we only sketch them. 
In Proposition \ref{lemsmallos}, we established, under $\phni \big( \cdot \big\vert \omega(T) > u \big)$, that $O_i$ has probability at least $c$ of being of size at most $c^{-1}$, where a small $c>0$ may be chosen uniformly in $u>0$ and the $i$-absent foliage. Given this, we arrive at the conclusion that $O_i = \emptyset$ has positive probability under the same law, by
 considering the event that the optional foundation be a finite length longer than it is in the case that $\vert V(O_i) \vert \leq c^{-1}$. This requirement on the optional foundation has a bounded probabilistic cost, but it is enough to ensure that $O_i$ may be empty with the condition $\omega(T) > u$ being satisfied.

In the case that $h_1 = 0$, it is impossible that $O_i = \emptyset$. However, we may have $D(O_i) = 1$, and the preceding argument works to show that the probability of this outcome under $\phni$ is bounded away from zero, uniformly in the $i$-absent foliage. $\Box$\\
The following lemma is the tool that we will use to  obtain Proposition \ref{lemdivexp} from Proposition \ref{lemsmallos} and Lemma \ref{lememptyos}.
\begin{lemma}\label{expdec}
Let $\big\{ X_i: i \in \mathbb{N} \big\}$ and $\big\{ R_i: i \in \N \big\}$ 
denote two sequences of $\mathbb{N}$-valued random variables such that 
$X_i \leq R_i$ almost surely for each $i \in \N$, and such that there exists $c > 0$ for which, for each $n,k \in \mathbb{N}$,
\begin{equation}\label{xnk}
\mathbb{P} \Big( X_n \geq k  \Big\vert  \big\{ X_j :  0 \leq j < n \big\} \Big) \leq \exp \big\{ - ck \big\},
\end{equation}
$\sigma \big\{ X_j: j < n\big\}$-a.s.,
\begin{equation}\label{xizero}
\mathbb{P} \Big( X_n \in \{ 0 , 1 \} \Big\vert \big\{ X_j: 0 \leq j < n \big\} \Big) \geq c,
\end{equation}
$\sigma \big\{ X_j: j < n\big\}$-a.s, and
\begin{equation}\label{rikineq}
\P \Big( R_n \geq k \Big\vert \big\{ R_j : 0 \leq j < n \big\}  \Big) \leq \exp \big\{ - ck \big\},
\end{equation}
$\sigma \big\{ R_j: j < n \big\}$-a.s.
 
Set $D_n = \big\{ n + 1 , \ldots, n + X_n - 1 \big\}$ if $X_n \geq 2$, and $D_n = \emptyset$ otherwise. 
Set $D = \bigcup_{n=0}^{\infty} D_n$, and $Y = \inf \big\{ i \geq 1: i \in
D^c \big\}$. (Note that $0 \in D^c$, so that $Y$ denotes the second smallest
member of $D^c \cap \N$.)
Further, set
$$
Z = \sum_{i=0}^{Y-1} R_i.
$$
Then there exists $\hat{c} > 0$, which is determined by $c$, such that, for each $k \in \mathbb{N}$,
\begin{equation}\label{eqzbound}
\mathbb{P}  \Big( Z \geq k  \Big) \leq \exp \big\{ - \hat{c} k \big\}.
\end{equation}
\end{lemma}
\noindent{\bf Proof of Proposition \ref{lemdivexp}.} 
Note firstly that Corollary \ref{cordecom} implies that it suffices to prove the statement with the choice $i=1$.

Let $T$ be a sample of the measure $\P_{h,\nu,u}$.
Recall from Definition \ref{defdecomplem} that, in the $\fso$-decomposition of $T$, the offshoot $O_i$ has root $\chi_i$. 
With a view to applying Lemma \ref{expdec} with the choice $\P = \P_{h,\nu,u}$, we set $X_i = D(O_i)$ and $R_i = \big\vert V(O_i) \big\vert$.

Deferring for a
moment the verification that the hypotheses of Lemma \ref{expdec} are
satisfied for this choice of sequences $\big\{ X_i: i \geq 0 \big\}$
and  $\big\{ R_i: i \geq 0 \big\}$, we note the reason that the lemma is
applicable: consider the set $D^c$ 
appearing in Lemma \ref{expdec}. Then the intersection of the vertex-set $V(S)$ of the spine with the set of cutpoints of $T$
consists precisely of those $\chi_i$, $0 \leq i \leq \vert S \vert$, for which $i \in D^c$. 
As such, the quantity $Y$ in Lemma \ref{expdec} is the index $i$ of
the lowest indexed cutpoint $\chi_i$ in $V(S) \setminus \{ \chi_0 \}$. 
Let $C^*$ denote the (renewal-decomposition) component  of $T$ whose root is $\chi_0$.
From the above, we see that
$$
 V\big( C^* \big) \subseteq \bigcup_{i=0}^{Y-1} V \big( O_i \big) \, \cup \,
 \big\{ \chi_Y \big\},
$$
whence
\begin{equation}\label{vcry}
\big\vert V\big( C^* \big) \big\vert \leq \sum_{i=0}^{Y-1} R_i \, + \, 1,
\end{equation}
%The constant $\hat{c}$ in the conclusion of Lemma \ref{expdec} depends only
%on the constant $c > 0$ appearing in the statement of the result. We
We find from Lemma \ref{expdec} and (\ref{vcry}) that
\begin{equation}\label{eqchatc}
\P_{h,\nu,u} \Big( \big\vert V(C^*) \big\vert \geq k \Big) \leq \exp \big\{
- \hat{c} k \big\},
\end{equation}
for all $u > 0$. Note that, if the foundation $F$ of $T$ has at least one edge, then, necessarily, $C_1$ consists of a single edge (in the case that $h_1 > 0$). If $F = \emptyset$, then $C_1 = C^*$. Hence,
 for $k \geq 2$, and for all $u > 0$, 
$$
\P_{h,\nu,u} \Big( \big\vert V(C_1) \big\vert \geq k \Big) \leq 
\P_{h,\nu,u} \Big( \big\vert V(C^*) \big\vert \geq k, 1\!\!1_{F = \emptyset} \Big) \leq 
\exp \big\{
- \hat{c} k \big\}.
$$
(A trivial modification is needed in the case that $h_1 = 0$.) This completes the proof, subject to checking
that
the hypotheses of Lemma \ref{expdec} are satisfied for the present choice of sequences and for $\P = \P_{h,\nu,u}$. In checking this, we must ensure that
we do so for a choice of the constant $c > 0$ in the hypotheses that is valid for
all $u > 0$, since we are claiming that there exists $\hat{c} > 0$ such that (\ref{eqchatc}) holds for all $u > 0$. 

Note that, for each of (\ref{xnk}), (\ref{xizero}) and (\ref{rikineq}), the event on which we condition is determined by the data in an $n$-absent foliage. Note also that the bound $X_i \leq R_i$, for any $i \geq 0$, is trivial. As such, 
 (\ref{xnk})
 and (\ref{rikineq}) are implied by Proposition \ref{lemsmallos}, and
(\ref{xizero}) by Lemma \ref{lememptyos}. $\Box$ \\
\noindent{\bf Proof of Lemma \ref{expdec}.}
We firstly show that, for some $c' > 0$,
\begin{equation}\label{yk}
\P \Big( Y \geq k \Big) \leq \exp \big\{ - c' k \big\}.
\end{equation}
To derive (\ref{yk}), it suffices 
to assume that $\big\{ X_i : i \in \mathbb{N} \big\}$
are independent and identically distributed, with $X_1$ having the law
$\chi = \chi_{c_0}$ for some $c_0 \in (0,1)$, that is determined by $c > 0$, where
$\chi$ on $\big\{ 1,\ldots,\infty \big\}$
has the form 
$$
\chi(i) = \big( 1 - c_0 \big) c_0^{i-1} 
$$
for $i \geq 1$. This is because we may choose $c_0 \in (0,1)$ 
such that, for each $n \in \mathbb{N}$,
the conditional distribution of $X_n$ given $\big\{ X_i : i < n \big\}$
is stochastically dominated by $\chi$, $\sigma \big\{ X_i : i < n
\big\}$-a.s.

Note that there exists $c_1 > 0$, determined by $c$, such that, for each $i \geq 1$,
\begin{equation}\label{idc}
\P \Big(  i \in D^c \Big) \geq c_1.
\end{equation}
Indeed, 
by use of (\ref{xnk}), there exists $K \in \N$ such that, for all $i \in
\N$,
$$
\P \Big( \bigcup_{j = K+1}^{i}  \big\{ X_{i-j} > j \big\}  \Big)
 \leq \frac{1}{2}.
$$
Conditional on the occurrence of $\bigcap_{j=K+1}^{i} \big\{
X_{i-j} \leq j\big\}$,
there is, by (\ref{xizero}), probability at least $c^K$ of
$\bigcap_{j=1}^K \big\{ X_{i-j} \leq j \big\}$.
Noting that 
$$
\Big\{ i \in D^c \Big\} = \bigcap_{j=1}^{i} \Big\{ X_{i-j} \leq j \Big\},
$$
we see that, as required for (\ref{idc}),
$$
\P \Big( i \in D^c \Big) \geq \frac{1}{2} c^K
$$
for such a value of $K$.

To show (\ref{yk}), firstly define $y_1 = \sup D_0 + 1$ if $D_0 \not=
\emptyset$ and $y_1 = 1$ otherwise. This quantity acts
as a candidate for the smallest element of $D^c$ exceeding $0$.
Note that, for each $n \geq 1$,
\begin{eqnarray}
& & \P \Big( y_1 \in D^c \Big\vert y_1 = n+1 \Big)
 = \P \Big( \bigcap_{j=1}^n \big\{ X_{n + 1 - j} \leq j \big\} \Big\vert X_0 =
 n + 1 \Big) \nonumber \\
 & = & \P \Big( \bigcap_{j=1}^n \big\{ X_{n  - j} \leq j \big\} \Big) 
 = \P \Big( n \in D^c \Big) \geq c_1, \nonumber
\end{eqnarray}
the second equality by the assumption that $\big\{X_i: i \in \N \big\}$ are independent and identically distributed, the inequality by (\ref{idc}).

If $y_1 \in D$, we seek a second candidate $y_2$ for the second smallest positive element
of $D^c$. Let $z_1 \in \{ 1, \ldots, y_1 - 1  \}$
be minimal such that $y_1 \in D_{z_1}$, that is, such that
$X_{z_1} \geq y_1 - z_1 + 1$.
We set $y_2 = z_1 + X_{z_1}$ (which is $\sup D_{z_1} + 1$).
Note that, from the form of the law $\chi$,
$X_{z_1} - \big( y_1 - z_1  \big)$, conditionally on $y_1 \in D$, has the
distribution of $\chi$. However, conditionally on $y_1 \in D$, we have that
 $y_2 - y_1 = X_{z_1} - \big( y_1 - z_1 \big)$, so that $y_2 - y_1$ has law 
$\chi$.
Note that 
$$
y_2 \geq \sup D_i  + 1
$$  
for any $i \leq z_1$, 
while the conditional distribution of 
$\big\{ X_{z_1 + j} : j \geq 1 \big\}$ is the same as the unconditioned
one. 
Thus, for any $n \geq 2$,
$$
\P \Big( y_2 \in D^c \Big\vert y_1 \in D , y_2 - z_1 = n  \Big)
 = \P \Big( n-1 \in D^c \Big) \geq c_1,
$$
the inequality by (\ref{idc}). Note also that $y_2 > y_1 > z_1$,
so that $y_2 - z_1 \geq 2$, implying that
 $\P \big( y_2 \in D^c \big\vert y_1 \in D \big) \geq c_1$.

We may iterate this procedure, forming a sequence $\big\{ y_i: i \in
\N\big\}$ such that $\big\{ y_{i+1} - y_i  : i \in \N \big\}$
is a sequence of independent terms, each having law $\chi$, 
with  each $y_i$, $i \in \N$, having probability at least $c_1$ of
belonging to $D^c$, conditional on all previous terms belonging to $D$, and
on any specific values for these previous terms.
As such, we have found that $Y$ is stochastically dominated by a sum of
independent geometric random variables, the number of summands being an
independent geometric random variable, whence (\ref{yk}). Note further that
$c'$ is determined by $c_0$ and $c_1$, and, thus, by $c$.

Turning now to (\ref{eqzbound}), 
note that, for any $c_2 > 0$,
$$
 \Big\{ Z \geq k \Big\}
 \subseteq \Big\{ \sum_{i=0}^{c_2 k} R_i > k \Big\} \cup \Big\{ Y \geq c_2 k \Big\},
$$
so that
\begin{equation}\label{eqzbdcomp}
\P \big( Z \geq k \big) \leq \P \Big( \sum_{i=0}^{c_2 k} R_i > k \Big)
 + \exp \big\{ - c' c_2 k \big\},
\end{equation}
by (\ref{yk}). In estimating $\P \big( \sum_{i=0}^{c_2 k} R_i > k \big)$, we may assume, due to(\ref{rikineq}), 
that $\big\{ R_i: i \in
\N \big\}$ is an independent and identically distributed sequence of random
variables whose law has an exponentially decaying tail, with constant in
the exponential determined by $c$.
Choosing $c_2 < \frac{1}{\E R_1}$, we obtain
\begin{equation}\label{eqrik}
 \P \Big( \sum_{i=0}^{c_2 k} R_i > k \Big) \leq \exp \big\{ - c_3 c_2 k \big\}.
\end{equation}
by means of the elementary large deviations bound
$$
\P \Big( \sum_{i=1}^m R_i > \big( \E R_1 + \epsilon \big) n \Big)
 \leq \big\{ - \psi(\epsilon) m \big\},
$$
(with $\psi(\epsilon) > 0$ for $\epsilon > 0$). Applying (\ref{eqrik}) to (\ref{eqzbdcomp}), we obtain (\ref{eqzbound}). This completes the proof of Lemma \ref{expdec}. $\Box$
\end{subsection}
\end{section}
\bibliographystyle{plain}
\bibliography{trwbib}
\end{document}